\documentclass[onefignum,onetabnum,onealgnum,oneeqnum,reqno]{siamart250211}
\usepackage[papersize={7.5in,10in},margin=0.975in,footskip=0.3in]{geometry}

\usepackage[final]{microtype}
\usepackage{amssymb,amsmath}
\usepackage{xcolor}
\usepackage{enumitem}
\usepackage{orcidlink}
\usepackage{booktabs}
\usepackage{multirow}
\usepackage{bm}

\setlist[enumerate]{leftmargin=2em,label=(\roman*),topsep=0.5em,parsep=0.25em}
\setlist[itemize]{leftmargin=1.2em,topsep=0.25em,parsep=0.25em}

\definecolor{bluey}{HTML}{008fba}
\hypersetup{urlcolor=bluey,linkcolor=bluey,citecolor=bluey}

\usepackage{caption}
\DeclareCaptionFont{smashmath}{\lineskiplimit=-\maxdimen}
\definecolor{bluesection}{HTML}{1e4b5e}
\makeatletter
\renewcommand\section{\@startsection{section}{1}{.25in}{1.3ex \@plus .5ex \@minus .2ex}{-.5em \@plus -.1em}{\reset@font\normalsize\bfseries\color{bluesection}}}
\renewcommand\subsection{\@startsection{subsection}{2}{.25in}{1.3ex\@plus .5ex \@minus .2ex}{-.5em \@plus -.1em}{\reset@font\normalsize\bfseries\color{bluesection}}}
\renewcommand\subsubsection{\@startsection{subsubsection}{3}{.25in}{1.3ex\@plus .5ex \@minus .2ex}{-.5em \@plus -.1em}{\reset@font\normalsize\bfseries\color{bluesection}}}
\makeatother

\usepackage{algpseudocodex}
\AtBeginEnvironment{algorithmic}{\footnotesize}
\crefname{figure}{Fig.\!}{Figs.\!}
\Crefname{figure}{Fig.\!}{Figs.\!}
\crefname{section}{\S\!}{\S\!}
\crefname{subsection}{\S\!\!}{\S\!\!}
\newcommand{\R}{\mathbb{R}}
\newcommand{\D}{\mathbb{D}}
\newcommand{\Id}{I}
\DeclareMathOperator{\diag}{diag}
\DeclareMathOperator*{\argmin}{arg\,min}
\newcommand{\vu}{\mathbf u}
\newcommand{\vn}{\mathbf n}
\newcommand{\vf}{\mathbf f}
\newcommand{\vg}{\mathbf g}
\newcommand{\vh}{\mathbf h}
\newcommand{\trans}{{\mathsf{T}}}
\newcommand{\jump}[1]{[\![ #1 ]\!]}
\newcommand{\polydeg}{{\smash{\wp}}}
\newcommand{\trimoverbar}[1]{\mkern 1.5mu\overline{\mkern-1.5mu#1\mkern-1.5mu}\mkern 1.5mu}
\newcommand{\SM}{SM}
\newcommand{\SMx}{Supplementary Material}

\newcommand{\Jmode}{\includegraphics[height=1.7ex]{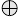}}
\newcommand{\GSmode}{\includegraphics[height=1.7ex]{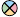}}
\newcommand{\FFmode}{\includegraphics[height=1.7ex]{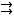}}
\newcommand{\RFmode}{\includegraphics[height=1.7ex]{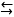}}
\newcommand{\FRmode}{\includegraphics[height=1.7ex]{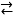}}
\newcommand{\RRmode}{\includegraphics[height=1.7ex]{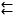}}

\newcommand{\getfig}[1]{\csname orb-page-#1\endcsname}
\makeatletter
\@namedef{orb-page-T7.Nx.P1}{1}
\@namedef{orb-page-T0.Nx.P1}{2}
\@namedef{orb-page-T4.Nx.P2}{3}
\@namedef{orb-page-T4.Nx.P3}{4}
\@namedef{orb-page-T0.Nx.P2.orbA}{5}
\@namedef{orb-page-T0.Nx.P3.orbA}{6}
\@namedef{orb-page-T0.Nx.P4.orbA}{7}
\@namedef{orb-page-T0.Nx.P2.orbB}{8}
\@namedef{orb-page-T0.Nx.P3.orbB}{9}
\@namedef{orb-page-T0.Nx.P4.orbB}{10}
\@namedef{orb-page-T0.Nx.P2.orbC-}{11}
\@namedef{orb-page-T0.Nx.P2.orbC+}{12}
\@namedef{orb-page-T0.Nx.P3.orbC-}{13}
\@namedef{orb-page-T0.Nx.P3.orbC+}{14}
\@namedef{orb-page-T0.Nx.P4.orbC-}{15}
\@namedef{orb-page-T0.Nx.P4.orbC+}{16}
\@namedef{orb-page-T5a.Nx.P3}{17}
\@namedef{orb-page-T5b.Nx.P3}{18}
\@namedef{orb-page-T2.Nx.P3.orbAbcP}{19}
\@namedef{orb-page-T2.Nx.P4.orbAbcP}{20}
\@namedef{orb-page-T2.Nx.P3.orbAbcD}{21}
\@namedef{orb-page-T2.Nx.P3.orbAbcN}{22}
\@namedef{orb-page-T2.Nx.P3.orbD1}{23}
\@namedef{orb-page-T2.Nx.P3.orbD2}{24}
\@namedef{orb-page-T2.Nx.P3.orbB}{25}
\@namedef{orb-page-T2.Nx.P3.orbC-}{26}
\@namedef{orb-page-T2.Nx.P3.orbC+}{27}
\@namedef{orb-page-T2.Nx.P3.orbEwx}{28}
\@namedef{orb-page-T2.Nx.P3.orbEgx}{29}
\@namedef{orb-page-T2.Nx.P3.orbD2.CS}{30}
\@namedef{orb-page-T2.Nx.P3.orbAbcP.CS}{31}
\@namedef{orb-page-T2.Nx.P3.orbC-.CS}{32}
\@namedef{orb-page-T63.Nx.P1}{33}
\makeatother
\newcommand{\putfig}[1]{\includegraphics[page=\getfig{#1}]{database.pdf}}


\newif\ifIncludeMAIN
\newif\ifIncludeSM
\IncludeMAINtrue
\IncludeSMtrue

\title{Cascading Smoothers for Multigrid}
\author{Robert I.~Saye \orcidlink{0000-0001-7256-6941}}
\date{\today}
\newcommand{\paperdate}{November 19, 2026}

\begin{document}
\setlength{\floatsep}{5mm}
\setlength{\textfloatsep}{5mm}

\ifIncludeMAIN

\makeatletter
\begin{center}
{\bfseries\MakeUppercase\@title\par}%
\vskip2mm%
{\small\textsc{\@author}\par}%
\vskip1mm%
{\sffamily\footnotesize Lawrence Berkeley National Laboratory, Berkeley, California, USA $\mid$ \texttt{rsaye@lbl.gov} $\mid$ \paperdate \par}%
\vskip5mm%
\end{center}
\makeatother

\begin{abstract}
Multigrid methods are among the most effective frameworks for solving large-scale sparse systems. However, achieving their hallmark linear scaling and rapid convergence crucially depends on an effective smoother algorithm, whose design is often highly problem-dependent. This paper develops a new approach, referred to as \textit{cascading smoothers} due to their operation as an ordered sequence of single-step block-diagonal smoothers. Each level in the cascade is optimised to maximally damp the output of prior steps via a Frobenius norm minimisation of the corresponding error propagators. In particular, we develop an additive (resp., multiplicative) formulation analogous to Jacobi (resp., Gauss-Seidel). Applied within a standard multigrid V-cycle, we show they are remarkably effective across a wide array of problems, including finite difference, finite element, and discontinuous Galerkin discretisations applied to Poisson, elliptic interface, and Stokes systems as well as multiphase variants. In every case, cascading smoothers closely match or significantly outperform their optimally-damped classical counterparts, yet require no parameter tuning apart from a few discrete solver choices. Additionally, the approach is highly parallelisable and robust to geometric and operator complexities such as unstructured meshes and high-contrast coefficients.
\end{abstract}



\section{Introduction}

Multigrid solvers are a cornerstone of scientific computing, widely recognised as some of the most efficient and versatile across a wide array of large-scale problems. Spanning fields from fluid dynamics and structural mechanics to plasma physics, geodynamics, and quantum mechanics (see, e.g., \cite{WESSELING2001311,10.1038/nature23911,ADAMS20106208,doi:10.1137/20M1348364,10.1145/2807591.2807675,RevModPhys.72.1041}), their principal domain of application lies in solving various classes of discretised partial differential equations (PDEs); we focus here on linear problems, generically represented as $Ax = b$. Among these, Poisson and Stokes constitute two of the canonical model problems at the core of multigrid research. Poisson problems represent the classic foundation where multigrid was initially forged \cite{FEDORENKO19621092,Brandt,Hackbusch} and continues to evolve, while Stokes problems serve as an archetype for coupled saddle-point systems \cite{BenziGolubLiesen2005,mgbook1}. %

Conceptually, multigrid algorithms use a hierarchical, multiscale error smoothing process that exploits the spectral properties of relaxation schemes. Classical methods, such as Jacobi and Gauss-Seidel, rapidly damp high-frequency, oscillatory errors but are ineffective on low-frequency, smooth errors. Multigrid resolves this issue by transferring the smoothed error (via the residual) onto a hierarchy of progressively coarser spaces. On coarser scales, the error effectively appears higher frequency and can be efficiently damped via relaxation, thereby forming a correction that is interpolated back to finer scales. Multigrid algorithms rely on three interlinked elements: a smoother method executed before and after coarse-grid corrections, intergrid transfer operators for residual restriction and correction interpolation, and a sufficiently consistent coarse-scale representation of $Ax = b$. When these components work in concert, multigrid can achieve optimal solver complexity, i.e., computational cost scaling linearly with the number of unknowns \cite{mgtutorial,mgbook1,mgbook3}. %

The two main pillars of multigrid design are geometric multigrid (GMG) and algebraic multigrid (AMG). Rather than relying on physical geometry, AMG frameworks operate solely on the entries of $A$ to automatically construct the hierarchy, transfer operators, and coarse-grid equations. This process centres on discovering a coarse space capable of representing the algebraically smooth error components relative to a prescribed relaxation method. As such, AMG yields a problem-independent approach, however this generality is paired with a computationally expensive setup phase; for a review, see \cite{STUBEN2001281,mgbook1,XuZikatanov2017}. In contrast, GMG assumes the fine-grid problem arises from a mesh that can be hierarchically coarsened and leverages naturally defined intergrid transfer operators, e.g., piecewise polynomial interpolation in finite element methods. Accordingly, GMG is primarily concerned with the design of an effective smoother. While AMG offers broader generality, GMG is typically preferred whenever feasible owing to its minimal computational overhead \cite{doi:10.1137/1.9780898718133.ch10,STUBEN2001281,mgtutorial}. To advance this paradigm, we focus here on GMG and develop a new kind of smoother.

Given an input approximation $x_i$ and right-hand side $b$, a single step of a smoother yields $x_{i+1} = x_i - B (A x_i - b)$ where $B$ is a given linear operator; the corresponding error propagator is given by $I - BA$. 
The landscape of multigrid smoothers is exceptionally rich. The following non-exhaustive overview highlights the methods most relevant to our work, focusing on local schemes with point- or block-wise updates at the scale of individual mesh elements.
\begin{itemize}
    \item \textit{Jacobi and Gauss-Seidel \cite{doi:10.1137/24M1680428}.} These classical iterative methods are the prototypical multigrid smoothers and employ a damping or overrelaxation parameter $\omega$ for maximum effectiveness. Damped Jacobi corresponds to a scaled inverted diagonal, yielding the error propagator $I - \omega D^{-1} A$ where $D = \diag(A)$; damped Gauss-Seidel (also known as successive over-relaxation) operates similarly, except updates are overwritten in-place according to a predefined component ordering. Optimal damping parameters are intrinsically problem-specific (see, e.g., \cref{fig:omega}), though local Fourier analysis (LFA) can be highly effective at providing sharp estimates in some idealised settings \cite{mgbook1,mgbook2,doi:10.1137/19M1308669}. In practice, these classical smoothers remain exceptionally effective across many canonical multigrid applications.
    
    \item \textit{Scheduled Relaxation \cite{YANG2014695,YANG2017135}.} %
    Rather than fixing $\omega$ over multiple iterations, scheduled relaxation uses a sequence of step-dependent damping parameters $\omega_i$; %
    a typical schedule leverages a combination of under- and over-relaxation steps to maximally damp high-frequency error. %
    Runge-Kutta smoothers and related pseudo-time-stepping relaxation methods \cite{10.2514/6.1989-1933,birken2016study} employ a similar multi-stage approach, optimising stage coefficients over the complex spectra typical of advection-dominated regimes.

    \item \textit{Polynomial and Chebyshev Smoothers \cite{AdamsBrezinaHuTuminaro,doi:10.1137/100798806,10.1002/nla.2518}.} Coupled to a base preconditioner $M$ (such as Jacobi), polynomial smoothers compute a multi-stage update with error propagator of the form $p_m(MA)$, where $p_m$ is a degree-$m$ polynomial satisfying $p_m(0) = 1$. Chebyshev smoothers minimise the maximum magnitude of $p_m$ over a high-frequency spectral interval $[\lambda^\star, \lambda_{\max}]$, often estimated by Lanczos or power iteration; Lottes \cite{10.1002/nla.2518} recently showed that fourth-kind Chebyshev smoothers yield quasi-optimality without requiring ad-hoc estimates for $\lambda^\star$. Polynomial smoothers feature exceptional parallel efficiency. %
    However, their smoothing ability still depends on accurate spectral estimates, and Chebyshev is usually constrained to problems with spectra in a real-valued interval. %

    \item \textit{Sparse Approximate Inverse (SPAI) Smoothers \cite{1997doi:10.1137/S1064827594276552,2000doi:10.1137/S0895479899339342,2001doi:10.1137/S1064827500380623,2002BROKER200261}.} SPAI methods construct a smoother $B$ by minimising the Frobenius norm $\|I - BA\|_F$ subject to a prescribed sparsity structure on $B$. Static methods predefine a fixed pattern, typically matching that of $A$ or $A^2$ or their block-form counterparts, or sparsified versions thereof; dynamic methods start with a minimal sparsity pattern and adaptively refine it, e.g., by prioritising strongly-coupled degrees of freedom. SPAI methods are compelling for their highly scalable implementation, not least because every row of $B$ can be computed independently and concurrently; furthermore, their robustness to unstructured meshes, complex geometry, and discontinuous PDE coefficients makes them among the most black-box smoothers available. However, SPAI lacks a simple mechanism for targeting high-frequency error components \cite{2010doi.org/10.1155/2010/930218} and a minimalist sparsity structure may not always suffice for robust multigrid performance \cite{2001doi:10.1137/S1064827500380623,inferno}.

    \item \textit{Overlapping Schwarz and Patch Smoothers \cite{doi:10.1137/1034116,Pav93,JK11,Cas97}.} Grounded in domain decomposition and subspace correction theory \cite{doi:10.1137/1034116}, these methods solve local problems on overlapping patches such as vertex-star or extended element neighbourhoods. Patch updates are combined additively, via weighted averaging across shared degrees of freedom, or multiplicatively, often via coloured sweeps. Notably, for high-order discretisations, overlapping Schwarz smoothers can yield provable $p$-robustness (i.e., convergence rates bounded independently of polynomial degree), with recent developments targeting modern architectures, reduced complexity, and matrix-free formulations by exploiting sparse bases, tensorisation, and fast diagonalisation \cite{10.1145/3445791,LF05,BF22,CGKS24,Wic25a}. For saddle-point Stokes problems, Vanka-type smoothers \cite{Vanka1986,Vanka1986b,10.1002/nla.2306} can be highly effective, essentially functioning as overlapping-block Jacobi or Gauss-Seidel methods applied to small cell-local patches grouping pressure unknowns with their connected velocity degrees of freedom. %

\end{itemize}

Adding to this landscape, we develop here a new class of methods, referred to as \textit{cascading smoothers} in light of their construction and operation. In essence, cascading smoothers apply an ordered sequence of single-step block-diagonal smoothers; each smoother is computed by minimising the Frobenius norm of its corresponding error propagator, conditioned on the output of prior smoothers in the sequence. %
The innermost smoother damps the error modes it identifies as most energetic, potentially amplifying others; the second level maximally damps the output of the first; and so on through the cascade. In one sense, this approach mirrors that of a scheduled damping block-Jacobi method, except there are no damping parameters to tune and the block components are not constrained in the same way. In another sense, it mirrors a polynomial smoother as a multi-stage process involving successive applications of $A$, except it requires no explicit spectral-bound estimates or polynomial design (nor does the overall error propagator take polynomial form). In a third sense, cascading smoothers share the Frobenius-norm minimisation approach of SPAI, except they enforce a minimal block-diagonal sparsity pattern at each level of the cascade; consequently, each successive level refocuses the optimisation toward high-frequency error modes. Finally, overlapping-block cascading smoothers closely resemble overlapping Schwarz and patch smoothers, except local patch solves are replaced by Frobenius-optimal block-diagonal operators conditioned on prior residuals; accordingly, the cascading approach implicitly accounts for inter-patch coupling and enables, e.g., effective multigrid Stokes solvers with smaller blocks than classical Vanka methods. On the other hand, cascading smoothers incur construction costs that grow with block size and residual fill-in at larger depths, a trade-off we revisit later through residual truncation and possible sketching approaches.

We show that multigrid solvers using cascading smoothers are remarkably effective across a broad range of problems, requiring no parameter tuning apart from a few discrete choices, e.g., regarding their forward and reverse application. The paper is structured as follows. In \cref{sec:cs}, we develop cascading smoothers, including additive and multicoloured formulations, describe a simple prescaling method that balances the system prior to Frobenius norm minimisation, and outline their application in a standard multigrid V-cycle. Extensive numerical results are presented in \cref{sec:numerics}, evaluating multigrid performance across a variety of test cases including finite difference, finite element, and discontinuous Galerkin frameworks applied to Poisson, steady-state and unsteady Stokes problems, alongside high-contrast multiphase variants. A concluding discussion is given in \cref{sec:conclusion}. The \SMx{} (\SM{}) contains additional discussion as well as an expanded set of numerical results. %

\section{Cascading Smoothers}
\label{sec:cs}

Cascading smoothers apply an ordered sequence of single-step block-diagonal smoothers. 
For now, we consider non-overlapping block methods where each block maps to the degrees of freedom on each mesh element; e.g., discontinuous Galerkin methods have a natural blocking corresponding to the collective set of nodal/modal coefficients on each element, while finite difference Poisson problems might simply use one degree of freedom per block. We partition $Ax = b$ into blocks so that $\sum_j A_{ij} x_j = b_i$, where $A_{ij}$ denotes the $(i,j)$th block of $A$; typically, each block is square, but may, e.g., be rectangular in a mixed-degree finite element method.

Given an input approximation $x_0$ and a right-hand side $b$, a single application of a \textit{depth $\nu$ (additive) cascading smoother} is defined by the output $x_\nu$ of the sequence
\begin{equation} \label{eq:fwda} x_\ell := x_{\ell-1} - \Lambda_\ell (A x_{\ell-1} - b), \quad \ell = 1, \ldots, \nu, \end{equation}
where $(\Lambda_\ell)_{\ell=1}^\nu$ are given block-diagonal matrices. By way of analogy, $\nu$ steps of block-Jacobi is equivalent to $\cref{eq:fwda}$ except that every $\Lambda_\ell$ is replaced by the inverse of the block-diagonal part of $A$. In general, each $\Lambda_\ell$ will be different, and \cref{eq:fwda} defines the \textit{forward application} of the corresponding cascading smoother. The essence of this work is to show that a simple least-squares based construction of $\Lambda_\ell$ leads to remarkably effective multigrid solvers. To that end, we first consider a Jacobi-style formulation, followed by a multicoloured Gauss-Seidel-style variant; these are referred to as \textit{additive} and \textit{multiplicative} formulations, analogous to Schwarz domain decomposition methods. 

\subsection{Additive Formulation}

Suppose for the time being that we have appropriately prescaled the operator $A$. The additive formulation of a cascading smoother is then built through a sequence of simple least squares minimisation problems, starting with the innermost level-one smoother, $\Lambda_1 := \argmin_{\Lambda \in \D} \| \Id - \Lambda A \|_F$, afterwhich we solve for the level-two smoother, $\Lambda_2 := \argmin_{\Lambda \in \D} \| (\Id - \Lambda A) (\Id - \Lambda_1 A) \|_F$, and then the level-three smoother, $\Lambda_3 := \argmin_{\Lambda \in \D} \| (\Id - \Lambda A) (\Id - \Lambda_2 A) (\Id - \Lambda_1 A) \|_F$, and so forth. Here, $\D$ denotes the set of block-diagonal matrices matching the block structure of $A$, $\Id$ is the identity matrix, and $\| \cdot \|_F$ denotes the Frobenius norm. In this sense, $\Lambda_\ell$ is the best possible block-diagonal smoother that minimises the Frobenius norm of its corresponding error propagator, $E_\ell := (\Id - \Lambda_\ell A) E_{\ell-1}$, conditioned on the error propagator of the prior smoothers in the sequence, $E_{\ell -1}$, with a base case of $E_0 := \Id$. (In the following, $E_\ell$ is also referred to as the \textit{residual} at the $\ell$-th level.) Note also that solving for $\Lambda_\ell$ is trivially parallelisable across its blocks: the Frobenius norm minimisation problem naturally decouples into independent subproblems, each taking the form of a standard overdetermined least squares problem.

\subsection{Multiplicative Formulation}

Analogous to a multicoloured Gauss-Seidel method, we can also design multicoloured cascading smoothers. Suppose we have $\chi$ different colours, let $\chi_k$ denote the set of elements with colour $k \in \mathbb{N} \cap [1,\chi]$, and let $\D_k \subseteq \D$ denote the set of block-diagonal matrices such that $(\Lambda \in \D_k) \Rightarrow (\Lambda_{ii} = 0\ \forall\ i \notin \chi_k)$. Then, we cycle through the colours (in some predefined order) and recursively solve for $\Lambda_{\ell,k} := \argmin_{\Lambda \in \D_k} \| (\Id - \Lambda A) E_{\ell,k-1} \|_F$, where $E_{\ell,k} := (\Id - \Lambda_{\ell,k} A) E_{\ell,k-1}$ and $E_{\ell,0} := E_{\ell-1,\chi}$ if $\ell > 1$ else $E_{1,0} := \Id$. Given an input approximation $x_0$ and a right-hand side $b$, a single (forward) application of a multiplicative cascading smoother is defined by the intuitive analogue of \cref{eq:fwda}, i.e., by the output $x_{\nu,\chi}$ of
\begin{equation} \label{eq:fwdm} x_{\ell,k} := x_{\ell,k-1} - \Lambda_{\ell,k} (A x_{\ell,k-1} - b), \quad \text{where} \quad x_{\ell,0} := \begin{cases} x_{\ell - 1, \chi} & \text{if } \ell > 1, \\ x_0 & \text{if } \ell = 1. \end{cases} \end{equation}
In this sense, $\Lambda_{\ell,k}$ corresponds to a block-diagonal smoother that acts only on the elements of colour $k$, while preserving the state on all other elements. (For example, $\nu$ steps of multicoloured block Gauss-Seidel is equivalent to \cref{eq:fwdm} except that every $\Lambda_{\ell,k}$ is replaced by the inverse of the block-diagonal part of $A$ restricted to colour $k$.) Similar to the additive formulation, solving for $\Lambda_{\ell,k}$ is trivially parallelisable across all of its nontrivial diagonal blocks. None of this construction actually depends on how the elements/blocks are coloured; in fact, the formulation collapses to an additive cascading smoother in the degenerate case that all elements have the same colour. However, intuitively we may expect that colouring in the same way as conventional multicolour Gauss-Seidel will provide the most benefits, e.g., in terms of convergence rates, parallelisation, construction cost, etc.; we discuss this more in later sections.

\subsection{Prescaling}
\label{sec:prescale}

In many applications, the entries of $A$ can vary over multiple orders of magnitude, e.g., in high-contrast multiphase problems, in Stokes problems that contain multiscaled differential operators, on highly unstructured mesh geometry, etc. If left untreated, this can skew the above uniformly-weighted least squares minimisation problems and result in ineffective multigrid smoothers. Fortunately, these aspects can be treated via a simple symmetric diagonal prescaling: instead of applying the above construction directly to $A$, we instead apply it to $W\!AW$, where $W$ is an appropriately-chosen diagonal matrix; once computed, the smoothers are then transformed back to the unscaled system by replacing each $\Lambda$ with $W\!\Lambda W$.\footnote{Corresponding to the problem $Ax = b$ we have the transformed problem $\tilde{A} \tilde{x} = \tilde{b}$, where $\tilde{A} = W\!AW$, $\tilde{x} = W^{-1} x$, and $\tilde{b} = W b$; applying a block-diagonal smoother $\tilde{\Lambda}$ to the latter is equivalent to applying the block-diagonal smoother $W\tilde{\Lambda}W$ to the former.} We have used two simple approaches to define the scaling $W$, depending on whether the problem is Poisson-like or Stokes-like.

\subsubsection{Elliptic problems} On scalar-valued elliptic problems, $A$ is symmetric positive (semi)definite and guaranteed to have positive diagonal entries. We prescale $A$ so that it has unit diagonal entries, i.e., by applying $W := \diag(A_{ii}^{\smash{-1/2}})$.

\subsubsection{Stokes problems} \label{sec:prescaleStokes} Stokes problems take the form
\[ Ax = b \quad \Leftrightarrow \quad \begin{pmatrix} \mathcal{M}_\rho + \mathcal{A}_\mu & \mathcal{G} \\ \mathcal{D} & \mathcal{P} \end{pmatrix} \begin{pmatrix} \vu_h \\ p_h \end{pmatrix} = \begin{pmatrix} b_{\vu} \\ b_{p} \end{pmatrix}\!, \]
where $x = (\vu_h, p_h)$ contains the velocity and pressure degrees of freedom, $b_\vu$ supplies the source data for the Stokes momentum equations, and $b_p$ supplies the source data for the divergence constraint;\footnote{In many applications of Stokes solvers, a divergence-free velocity field is sought, i.e., $b_p$ is identically zero on the finest mesh. Nevertheless, the multigrid solvers developed here allow for nonzero $b_p$ on every level of the hierarchy, so that we can smooth and perform coarse-grid corrections for the velocity and pressure error components, simultaneously.} meanwhile, $\mathcal{M}_\rho$ is a density-weighted mass matrix (zero for steady-state problems), $\mathcal{A}_\mu$ is a viscosity-weighted Laplacian-like viscous operator, %
$\mathcal{G}$ is a pressure gradient operator, and $\mathcal{D} = \mathcal{G}^\trans$ is a negative-divergence operator.\footnote{$\mathcal{P}$ is an operator that might be necessary for pressure stabilisation, depending on the specifics of the discretisation. For all applications presented in this work, no stabilisation is needed and $\mathcal{P} \equiv 0$ throughout. When $\mathcal{P} \neq 0$, its scaling is more complex and relates to the Schur complement of the Stokes system; the methods of \cref{sec:prescaleStokes} implicitly assume that the induced rescaling $\diag(w_p) \mathcal{P} \diag (w_p)$ is well balanced. While further exploration of $\mathcal{P} \neq 0$ is left for future work, we note that the approach proves effective in preliminary tests of stabilised equal-degree LDG Stokes formulations (in 2D and 3D across various polynomial degrees).} These operators (locally) scale as $\mathcal{M}_\rho \sim \rho h^d$, $\mathcal{A} \sim \mu h^{d-2}$, $\mathcal{G} \sim h^{d-1}$, and $\mathcal{D} \sim h^{d-1}$, where %
$h$ is the mesh element size. Note that the dependence on $h$ means that the length scales of the domain, or even the particular level of the multigrid hierarchy, can affect the relative strengths of these operators. To balance them, we define $W := \diag(w_\vu, w_p)$, as follows:
\begin{enumerate}
    \item First, we define $w_\vu$ so that the balanced form of $\mathcal{A} := \mathcal{M}_\rho + \mathcal{A}_\mu$ has unit diagonal weighted by the relative $L^1$ row norm; i.e., we choose $w_\vu$ so that the $i$th diagonal of $\diag(w_\vu) \mathcal{A} \diag(w_\vu)$ is equal to $\| \mathcal{A}_{i,:}\|_1 / \mathcal{A}_{ii}$, where $\mathcal{A}_{i,:}$ denotes the $i$th row of $\mathcal{A}$.
    \item Second, we define $w_p$ so that every row of the balanced form of $\mathcal{D}$ has unit max-norm, i.e., each row of $\diag(w_p) \mathcal{D} \diag(w_\vu)$ has maximum absolute value $1$. (Equivalently, since $\mathcal{D} = \mathcal{G}^\trans$, the prescaled form of $\mathcal{G}$ is such that it has unit max-norm columns.)
\end{enumerate}
This approach was chosen for its near-universal effectiveness across a multitude of Stokes problems, including steady-state, viscous-dominated and near-inviscid unsteady problems, as well as various discretisation frameworks %
including mixed-degree and equal-degree discontinuous Galerkin methods and Taylor-Hood finite element methods. There is one exception, corresponding to a finite difference MAC staggered grid discretisation; in this setting, it is more effective to normalise the rows (resp., columns) of $\mathcal{D}$ (resp., $\mathcal{G}$) by the $1$-norm instead of the max-norm. Additional discussion on this scaling approach is given in the \SM{}.

\subsection{Multigrid Schemes}
\label{sec:vcycle}

\begin{figure}[!t]%
\centering%
\begin{minipage}{101mm}%
\begin{algorithm}[H]
	\caption{Multigrid V-cycle $V({\mathcal E}_h, x_h, b_h)$ on mesh ${\mathcal E}_h$ of the hierarchy.}
	\begin{algorithmic}
		\If{${\mathcal E}_h$ is the bottom level}
			\State Solve $A_h x_h = b_h$ using the bottom level direct solver.
		\Else
			\State Apply pre-smoother to $x_h$.
			\State Compute restricted residual, $r_{2h} := (I_{2h}^h)^\trans (b_h - A_h x_h)$.
			\State Solve coarse grid problem, $x_{2h} := V({\mathcal E}_{2h}, 0, r_{2h})$.
			\State Interpolate and correct, $x_h \leftarrow x_h + I_{2h}^h x_{2h}$.
			\State Apply post-smoother to $x_h$.
		\EndIf
		\State{\textbf{return} $x_h$.}
	\end{algorithmic}%
	\label{algo:vcycle}%
\end{algorithm}%
\end{minipage}%
\end{figure}

The multigrid solvers we consider here are built using a standard geometric multigrid V-cycle. On the finest level we consider uniform Cartesian grids as well as quadtree/octree meshes for problems with adaptive mesh refinement; in each case, there is a naturally-defined nested mesh hierarchy $\mathcal{E}_h, \mathcal{E}_{2h}, \ldots$. The interpolation operator $I_{2h}^h$ and restriction operator $R_h^{2h} = (I_{2h}^h)^\trans$ depend on the specific test application, but are standard, and will be noted in context. All of our examples use coarse-grid operators that are equivalent to rediscretising on every level of the hierarchy; in many cases, this is equivalent to Galerkin coarsening, i.e., $A_{2h} = (I_{2h}^h)^\trans A_h I_{2h}^h$. Finally, a direct solver is used on the bottom level of the grid hierarchy, typically consisting of just one element spanning the entire domain. %
With this setup, the template of a standard V-cycle is shown in \cref{algo:vcycle}. After specifying a pre- and post-smoothing method, the invocation of the V-cycle on the finest mesh with an initial guess of zero yields a linear operator acting on the given right-hand side $b$; this linear operator is denoted in the following by $V$.

To use a cascading smoother in a multigrid method, we have some options regarding the ordering of the underlying single-step smoothers. In particular, \Cref{eq:fwda} and \cref{eq:fwdm} define the \textit{forward applications} of the additive and multiplicative version, respectively. We also consider their adjoint-counterpart \textit{reverse applications}, which inverts the order and transposes $\Lambda$. Explicitly, given an approximate solution $x_0$ to $Ax = b$, the reverse application of an additive cascading smoother $(\Lambda_\ell)_{\ell=1}^\nu$ is given by the output $x_\nu$ of the sequence
\begin{equation} \label{eq:fwdar} x_\ell := x_{\ell-1} - \Lambda_{\nu-\ell+1}^\trans (A x_{\ell-1} - b), \quad \ell = 1, \ldots, \nu. \end{equation}
The analogue for a multiplicative cascading smoother $(\Lambda_{\ell,k})_{\ell,k=1}^{\nu,\chi}$ is given by \cref{eq:fwdm}, except the colour ordering is also reversed and is given by the output $x_{\nu,\chi}$ of
\begin{equation} \label{eq:fwdmr} x_{\ell,k} := x_{\ell,k-1} - \Lambda_{\nu - \ell + 1,\chi - k + 1}^\trans (A x_{\ell,k-1} - b), \text{ where } x_{\ell,0} := \begin{cases} x_{\ell - 1, \chi} & \text{if } \ell > 1, \\ x_0 & \text{if } \ell = 1. \end{cases} \end{equation}

We consider here two main classes of multigrid solvers: (i) using the same additive cascading smoother, for both pre- and post-smoothing, denoted by the symbol $\Jmode$; and (ii) using the same multicoloured cascading smoother, for both pre- and post-smoothing, denoted by the symbol $\GSmode$. For each of these classes, we consider two combinations of ordering: (a) forward application for both pre- and post-smoothing, denoted by $\FFmode$; and (b) reverse application for the pre-smoother and forward application for the post-smoother, denoted by $\RFmode$. The choice $\RFmode$ leads to a symmetric V-cycle, i.e., $V = V^\trans$ (a proof is given in the \SM{}), whereas $\FFmode$ leads to an asymmetric $V$. (Other combinations of ordering have also been tested, but they are not included in the presented results.\footnote{Numerical experiments clearly and universally demonstrated that $\FFmode$ and $\RFmode$ are the best two approaches out of all four possible orderings, i.e., $\RRmode$ and $\FRmode$ yielded multigrid solvers whose convergence rates are either on-par with, or eclipsed by, $\FFmode$ and/or $\RFmode$.}) In total, we consider four solver types: $\Jmode\FFmode$, $\Jmode\RFmode$, $\GSmode\FFmode$, and $\GSmode\RFmode$. In general, putting aside various implementation trade-offs and all else being equal, the multicoloured methods yield faster convergence rates than the additive methods. Regarding the ordering, we will see that the optimal choice varies by application.

\subsection{Implementation}
\label{sec:impl}

\begin{figure}[!t]%
\centering%
\begin{minipage}{96mm}%
\begin{algorithm}[H]
	\caption{Building a depth $\nu$ additive cascading smoother $(\Lambda_\ell)_{\ell=1}^\nu$.}
	\begin{algorithmic}[1]
        \State Compute an appropriate diagonal prescaling matrix $W$ (see \cref{sec:prescale}).
        \State Prescale $A$ by replacing $A \leftarrow W\!A W$.
        \State Define $E_0 := \Id$.
        \For{level $\ell = 1, \ldots, \nu$}
            \For{every element $i$ (possibly in \textbf{parallel})}
                \State Define $\Lambda_{\ell,ii}$ as the least squares solution of
                \Statex \qquad $\Lambda_{\ell,ii} (A E_{\ell-1})_{i,:} = (E_{\ell -1})_{i,:}.$ \label{algo:solve}
                \State Compute the $i$th block row of the residual,
                \Statex \qquad $(E_{\ell})_{i,:} := (E_{\ell-1})_{i,:} - \Lambda_{\ell,ii} (A E_{\ell-1})_{i,:}.$
                \State Optionally, truncate the residual, e.g., by removing all but the largest $\alpha$ blocks of $(E_{\ell})_{i,:}$. \label{algo:trunc}
        \EndFor
        \EndFor
        \State For every level $\ell$, replace $\Lambda_\ell \leftarrow W\! \Lambda_\ell W$.
        \State Restore $A$ to its unscaled form, $A \leftarrow W^{-1}\!A W^{-1}$.
	\end{algorithmic}%
	\label{algo:csa}%
\end{algorithm}%
\end{minipage}%
\end{figure}

For every level of the multigrid hierarchy and its associated operator $A_h$, we build a corresponding cascading smoother for use in the V-cycle. \Cref{algo:csa} outlines the basic implementation for an additive cascading smoother; a multicoloured formulation is the same except for an additional loop over the set of colours; %
an implementation that re-uses/overwrites state variables is given in the \SM{}. The depth, $\nu$, is a user-defined parameter, analogous to the number of pre- and post-smoothing steps of a conventional smoother; in practical settings, $\nu = 2$ or $\nu = 3$ typically gives the best results, and there is usually no benefit in going beyond $\nu = 4$, as will be shown in the next section.

We have yet to discuss the computational cost of building a cascading smoother. Each block of $\Lambda_\ell$ requires the solution of an overdetermined least squares problem of size $(k_\ell n^2) \times n$, where $n$ is the size of each block\footnote{While this analysis assumes a uniform block size $n$ for simplicity, cascading smoothers inherently support variable block sizes.} and $k_\ell$ is the number of nonzero blocks in the corresponding block-row of $A E_{\ell-1}$; using standard solver methods (e.g., normal equations, QR) this amounts to a construction cost of $\mathcal{O}(k_\ell n^3)$ per diagonal block of $\Lambda_{\ell}$. In the additive case, the residual $E_\ell$ has the same block-sparsity pattern as the $\ell$th power of $A$; for a typical PDE-based problem, it follows that $k_\ell \in \mathcal{O}(\ell^d)$, where $d$ is the spatial dimension, and the per-block-row cost of forming $A E_{\ell - 1}$ (and $E_\ell$, if necessary) is $\mathcal{O}(k_\ell n^3)$. The analysis for a multicoloured formulation is more subtle, as it depends on the colouring method: the worst-case scenario results in a strictly sequential dependency chain so that $E_{\ell,\chi}$ is dense after just one level; the best-case scenario is analogous to standard multicoloured Gauss-Seidel, whereby we choose the fewest possible colours that maximise parallelism within each colour set.

Whether additive or multiplicative, the predominant factor in overall construction cost of deeper smoothers is that of fill-in of the residual $E_\ell$. For many applications and/or shallow cascading smoothers, this is not a problem, e.g., whenever the smoother can be precomputed offline; see also the \SM{} for an empirical evaluation of execution times and runtime breakdown on representative examples. Applying a truncation or dropoff method offers a simple way to speed up construction, as suggested on line \ref{algo:trunc} of \cref{algo:csa}; i.e., before proceeding to the next level, remove all but the largest $\alpha \in \mathbb{N}$ blocks of the $E_\ell$. Provided $\alpha$ is large enough, the resulting cascading smoother yields practically similar multigrid solver speeds as $\alpha = \infty$; often $\alpha \approx 5$ in 2D and $\alpha \approx 7$ in 3D suffices. Besides truncation methods, a compelling alternative is to use sketching approaches, which may open up more advanced computational paradigms such as matrix-free implementations (e.g., where the system operator $A$ is a black box), GPU acceleration, and high-order finite element methods using tensorised evaluation---a topic we revisit in the concluding remarks.

\subsection{Extension to Overlapping Blocks}
\label{sec:overlapping}

While the non-overlapping block framework established above is effective across a wide range of problems, even for various kinds of saddle-point Stokes problems, some applications benefit from an overlapping-block approach. Conceptually, both the construction and application are the same, except for the temporary formation and subsequent fusion of overlapping blocks at each level. (Since the upcoming numerical studies focus primarily on non-overlapping methods, the reader is welcome to skip this subsection without loss of continuity.)

Given a collection of overlapping blocks, let $\Theta$ denote the corresponding tall-and-narrow binary matrix that replicates shared degrees of freedom, mapping input vectors to an expanded space in which these blocks become disjoint; each row of $\Theta$ has exactly one $1$. Let the vector $\theta$ denote the associated multiplicity count, given by the column sums of $\Theta$, i.e., $\theta_i := \| \Theta_{:,i} \|_1$. A single application of a \textit{depth $\nu$ overlapping-block additive cascading smoother} $(\Lambda_\ell)_{\ell=1}^\nu$ is defined via the following analogue of \cref{eq:fwda}:
 \[ x_\ell := x_{\ell-1} - \diag(\theta)^{-1}\Theta^\trans \Lambda_\ell\, \Theta (A x_{\ell-1} - b), \quad \ell = 1, \ldots, \nu, \]
i.e., we compute the residual as normal, temporarily replicate degrees of freedom, apply the block-diagonal smoother $\Lambda_\ell \in \bar{\D}$, and then fuse multivalued nodes via weighted averaging. Here, $\bar{\D}$ denotes the set of block-diagonal matrices associated with the overlapping blocks in the expanded space. To build $(\Lambda_\ell)_{\ell=1}^\nu$, we apply a nearly analogous process as the non-overlapping case: assuming $A$ has been appropriately prescaled (see \cref{sec:prescale}), then
\[ \Lambda_\ell := \argmin_{\Lambda \in \bar{\D}} \|(\Theta - \Lambda \Theta A) E_{\ell-1} \|_F, \quad \text{where} \quad E_\ell := \diag(\theta)^{-1} \Theta^\trans \bigl(\Theta - \Lambda_\ell \Theta A\bigr) E_{\ell-1}, \]
with a base case of $E_0 := I$. In particular, note that $\Lambda_\ell$ solves a Frobenius-norm minimisation problem in the expanded space, prior to collapsing the smoother residual back to the original space. This allows us to maintain the benefits of independent, decoupled subproblems, parallelisable across the set of diagonal blocks. Meanwhile, an overlapping-block multiplicative cascading smoother requires no weighted averaging; its implementation is essentially identical to that of the non-overlapping case, except that the block structure shifts in tandem with the active block (for the finer details, see the \SM{}).

\section{Numerical Experiments}
\label{sec:numerics}

In this section, we examine the effectiveness of cascading smoother multigrid solvers across a variety of test problems in 2D and 3D. Key results are highlighted here; additional test configurations, including boundary condition types, polynomial degrees, and problem coefficients, are given in the \SM{}.

\subsection{Preliminaries}
\label{sec:prelim}

It is a commonplace and highly effective strategy to accelerate a multigrid solver by using it as a preconditioner in a Krylov method, such as Conjugate Gradient (CG) and Generalised Minimal Residual (GMRES) \cite{demmel,saad}. We focus here on GMRES, in part for its suitability to Stokes problems, but also for its flexibility with asymmetric preconditioners, allowing us to test different forward/reverse applications of cascading smoothers.\footnote{Viability of a cascading smoother V-cycle as preconditioner for CG is discussed in the \SM{}.} %
Accordingly, throughout this work, our main driver solver is that of multigrid-preconditioned GMRES, using the V-cycle of \cref{sec:vcycle} as left-preconditioner.

We measure and compare multigrid solvers through two key metrics, $\eta$ and $\mathcal{C}$, each a function of the empirically-determined convergence rate, $\rho$. For each solver configuration and grid size, we randomly generate a right-hand side vector $b$ and solve $Ax = b$ using a zero initial guess $x_0 = 0$. Each iteration $i$ yields an improved solution $x_i$ with corresponding residual norm $r_i$; typically, $r_i$ is a nearly-constant fraction of $r_{i-1}$ so that $r_i \approx \rho^i\,r_{0}$, where $\rho$ is the convergence rate.\footnote{The components of $b$ are randomly drawn from the uniform distribution on $[-1,1]$ so that the test problem is likely to contain modes for which convergence is the slowest. The residual norm is the same as that used by left-preconditioned GMRES, i.e., $r_i := \|Vb - V\!Ax_i\|_2$. On a log-linear plot, the data points $(i, r_i)$ are well approximated by a straight line of negative slope; simple linear regression is used to compute the corresponding convergence rate $\rho$. Numerical tests typically have 5--10 such data points, except when convergence is so fast that just a few data points are collected before reaching the limits of double-precision arithmetic.} Our primary metric used to assess multigrid solver speed is then defined by
\begin{equation} \label{eq:eta} \eta := \frac{\log 0.1}{\log \rho} . \end{equation}
Here, $\eta$ represents the expected number of iterations needed to reduce the residual by a factor of $10$; e.g., $k\eta$ represents the number of iterations to reduce it by $10^k$. One reason to prefer $\eta$ over $\rho$ is that the former is more intuitively useful; e.g., a two-fold reduction in $\eta$ is the same as a two-fold reduction in iteration count, and so, all else being equal, two times faster.

Our second metric weighs the benefits of improved convergence rates against the increased cost of applying deeper cascading smoothers, analogous to increasing the number of pre- and post-smoothing steps of traditional smoothers. To simplify the analysis, we assume the overall solver cost is dominated by that of block-sparse matrix-vector multiplications; accordingly, let one cost unit be defined as the cost of a single block matrix-vector multiply. Suppose the operator $A$ has a stencil size $m$, i.e., $m$ nonzero blocks per block-row. To leading order, a single application of a depth $\nu$ cascading smoother costs $\nu m$ units per element, the same as $\nu$ steps of block-Jacobi/Gauss-Seidel.\footnote{The precise quantity depends on some implementation choices, e.g., whether Jacobi-style methods evaluate the residual in one step, followed by application of the inverse block-diagonal in a second step, or if the net update is performed in one sweep thereby saving one block multiply, etc. For our purposes, the leading order $\nu m$ metric suffices.} Factoring in the cost of residual evaluation ($m$ units per element), restriction ($r$ units per element) and correction interpolation ($i$ units per element), and summing over every level of the hierarchy, the total cost of a V-cycle is therefore $c_d (2 \nu m + m + r + i)$ units per fine-grid element, where the dimension-dependent prefactor $c_d := (1 - 2^{-d})^{-1}$ reflects that element counts typically reduce by a factor $2^d$ each level. Meanwhile, each GMRES iteration itself requires one residual calculation, adding another $m$ units per element. Finally, to reduce the residual by a factor of $10^k$, we need $\eta k$ iterations in total. For simplicity, let us reasonably suppose that the cost of interpolation and restriction can be neglected, so that we may arrive at the following conclusion: to reduce the residual by a factor of $10^k$, a V-cycle preconditioned GMRES method has cost proportional to $k m\, \mathcal{C}_d(\nu, \eta)$, where
\begin{equation} \label{eq:cost} \mathcal{C}_d(\nu, \eta) := \eta \Bigl( 1 + \frac{2\nu + 1}{1 - 2^{-d}} \Bigr). \end{equation}
This formula reflects that, although increasing $\nu$ is expected to decrease $\eta$, the overall benefit is a mixed compromise. Naturally, this analysis neglects a multitude of implementation factors, including arithmetic density, %
vectorisation, memory locality and cache effects, synchronisation costs when comparing additive and multiplicative methods, distributed and shared memory parallelism, communication costs of deep V-cycles, global reduce operations, and so forth. Nevertheless, $\mathcal{C}_d(\nu, \eta)$ serves to provide a leading-order approximation of the overall solver efficiency so that we may roughly compare different smoothers of different depths.

\subsection{Optimally-Damped Jacobi and Gauss-Seidel}
\label{sec:JGS}

Throughout this work, we compare against the classical methods of block-Jacobi and block-Gauss-Seidel (and for some Stokes problems, their Vanka counterparts). In many cases, Jacobi and/or Gauss-Seidel are known to yield extremely efficient multigrid solvers, and thus serve as a good baseline for comparison. %

\begin{figure}%
\centering%
\includegraphics[width=5.125in]{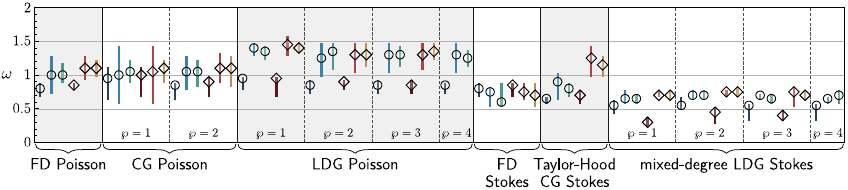}%
\caption[]{Optimal damping parameters for classical smoothers. Circles and blue/green shades correspond to 2D; diamonds and red/orange shades correspond to 3D. Within each group of three, the left-most point corresponds to $\mathsf{J}$, the middle point to $\smash{\mathsf{GS}\GSmode\RFmode}$, and right-most to $\smash{\mathsf{GS}\GSmode\FFmode}$ (or in the cases of FD Stokes and Taylor-Hood CG Stokes, their Vanka counterparts $\smash{\mathsf{V}\Jmode}$, $\smash{\mathsf{V}\GSmode\RFmode}$, and $\smash{\mathsf{V}\GSmode\FFmode}$, see \cref{sec:staggered} and \cref{sec:TH}); bars indicate ranges of $\omega$ where $\eta$ is within 10\% of optimal.}%
\label{fig:omega}%
\end{figure}

For each test application, we performed a parameter sweep to empirically determine an optimal damping/overrelaxation parameter $\omega$ for the following methods: (i) block-Jacobi, denoted by $\mathsf{J}$; (ii) multicoloured block-Gauss-Seidel, where the colour sweep ordering is kept the same for both the pre- and post-smoother, denoted by $\mathsf{GS}\GSmode\FFmode$; and (iii) multicoloured block-Gauss-Seidel, where the colour ordering is reversed for the pre-smoother, denoted by $\mathsf{GS}\GSmode\RFmode$. Both $\mathsf{J}$ and $\mathsf{GS}\GSmode\RFmode$ lead to symmetric V-cycle solvers, but not $\mathsf{GS}\GSmode\FFmode$. Their implementation is analogous to that of cascading smoothers, except that every block-diagonal matrix $\Lambda$ is replaced by the inverse block-diagonal part of $A$ and scaled by $\omega$. Using multigrid-preconditioned GMRES as the driver and a moderately-deep V-cycle (9--10 levels in 2D, 8 levels in 3D), the parameter sweep searches over the range $\omega \in (0,2)$ to find $\omega$ that minimises $\eta$;\footnote{Specifically, this calibration uses $\nu = 3$ pre- and post-smoothing steps as a conservative default across all applications, favouring cases where deeper smoothing is needed (e.g., those where $\nu = 1$ can fail). Although the optimal $\omega$ can mildly vary with $\nu$ (cf.~\cref{footnote:rbomega}), an exhaustive fine-tuning for every  $\nu$ is not our focus here; rather, the intent is to establish a meaningful, reasonably well-optimised classical smoother baseline.} the sweep also records a window of near-optimal values to assess sensitivity. The results are summarised in \cref{fig:omega}. It should be emphasised that there is not necessarily any pattern to discern in these results; rather, that each application (Poisson or Stokes; finite difference (FD), continuous Galerkin (CG), or discontinuous Galerkin (DG); 2D or 3D; polynomial degree; etc.) typically requires its own computationally-intensive parameter sweep that we would like to avoid in the broader context of designing multigrid solvers \cite{doi:10.1137/19M1308669}.

\subsection{Experimental Setup and Visualisation Guide}
\label{sec:guide}

For each test application, we measure multigrid performance across multiple grid sizes and every solver configuration: Jacobi and Gauss-Seidel smoothers, of varying pre- and post-smoother iteration counts $\nu$; additive and multiplicative cascading smoothers, of various depths $\nu$; and their various forward/reverse applications $\FFmode$ and $\RFmode$. To illustrate this wealth of data, we use visualisations of the form shown in \cref{fig:T7.Nx.P1}, \cref{fig:T4.Nx.P3}, etc. Each half corresponds to 2D and 3D; within each half:
\begin{itemize}
    \item The left-side plots $\eta$ and consists of two main rows, additive-style in the top row and multiplicative-style in the bottom. Divided primarily by $\nu \in \{1,2,3,4\}$, each panel is further subdivided into two sides: the left half corresponds to $\nu$ steps of Jacobi/Gauss-Seidel, whereas the right half corresponds to a depth $\nu$ cascading smoother. Within each panel, each individual sequence of markers plots the experimentally-determined multigrid convergence speed $\eta$ as measured on a grid of size $n \times n\,(\times n)$, one marker point per $n$, with $n = 4, 8, 16, 32, \ldots$, going from left-to-right. Sometimes a solver is slow to converge or may even be nonconvergent: whenever $\eta \geq 3$, we declare the method to be ineffective (at least, from a practical point of view) and indicate this by a faded marker placed outside the graph. Finally, open and closed circles distinguish the pre- and post-smoother ordering modes, as indicated in the legend; in the case of Jacobi, there is no ordering to specify, so there is only one sequence per $\mathsf{J}$ panel.
    \item Meanwhile, the right-side heat map visualises the $\nu$-weighted cost of each solver method, relative to the overall best. Each solver configuration yields a set of $\eta$ values from all the grid sizes tested; we use the third quartile, $\eta_{0.75}$, as a representative for the ``typical'' convergence rate among the slowest observed. Making use of \cref{eq:cost} to factor in the increased costs of larger $\nu$, the overall solver cost is defined by $\mathcal{C}_d(\nu, \eta_{0.75})$. We normalise these quantities by the minimum across all 28 configurations to define each solver's relative cost, $\mathcal{C}_\text{rel}$; e.g., $\mathcal{C}_\text{rel} \approx 2$ indicates a solver wall-clock time about twice that of the best. The heat map visualisation is shaded accordingly: best-performing solvers are highlighted via darker shading, non-competitive solvers are unshaded, while the blue and grey colouring is simply to help distinguish additive-style from multiplicative-style.
\end{itemize}

\subsection{Poisson Problems}

\subsubsection{Finite difference methods}
\label{sec:fd}

We begin with the canonical test case of a Poisson problem discretised by second-order finite differences on a uniform Cartesian grid. The corresponding operator $A$ is the standard 5-point ($d=2$) or 7-point ($d=3$) negative Laplacian with diagonal $2d/h^2$. We consider two approaches for the corresponding multigrid solver: (i) view the grid as node-centred and apply bilinear/trilinear interpolation in the V-cycle; alternatively, (ii) view the grid as cell-centred and apply piecewise constant interpolation, analogous to a finite volume method. We present results using the first approach here, and leave the second approach for the \SM{}.

\begin{table}%
\centering%
\sffamily\footnotesize%
\newcommand{\ssize}[2]{$#1 \times #2$}%
\begin{tabular}{lc|cccccc}
Cascading smoother formulation & $d$ & $\ell = 1$ & $\ell = 2$ & $\ell = 3$ & $\ell = 4$ \\
\midrule
\multirow{2}{*}{Additive, $\lambda_\ell$}
& 2D & 0.80 & 0.86 & 0.91 & 0.92 \\
& 3D & 0.86 & 0.88 & 0.90 & 0.93 \\
\midrule
\multirow{2}{*}{Multiplicative, $(\lambda_{\ell,\text{R}}, \lambda_{\ell,\text{B}})$}
& 2D & (0.80, 1.14) & (1.33, 1.45) & (1.55, 1.60) & (1.65, 1.68)\\
& 3D & (0.86, 1.12) & (1.24, 1.40) & (1.46, 1.53) & (1.59, 1.61) \\
\midrule
\end{tabular}%
\caption{Effective damping/overrelaxation coefficients for the cascading smoothers arising from the finite difference Poisson problem considered in \cref{sec:fd}.}%
\label{tab:autodamp}
\end{table}

We apply the trivial blocking of one grid point per block, and for multicoloured smoothers, the canonical red-black colouring. With periodic boundary conditions, every row of $A$ has the same stencil; it follows that every (purely) diagonal matrix of an additive cascading smoother $(\Lambda_\ell)_{\ell=1}^\nu$ collapses to a scaled identity matrix, $\Lambda_\ell = \lambda_\ell \diag(A)^{-1} = \lambda_\ell\, (2d/h^2)^{-1} I$, where $\lambda_\ell$ is constant depending only on $\ell$ and $d$, but not the cell size $h$. As such, we may reinterpret the additive cascading smoother as a damped Jacobi method, except with automatically-determined level-dependent damping coefficients $\lambda_\ell$. An analogous observation holds for multicoloured cascading smoothers: each level has two (purely) diagonal smoothers for the two colours, whose application is functionally the same as a damped multicoloured Gauss-Seidel method, except that each colour sweep has its own level-dependent, automatically-determined damping coefficient $\lambda_{\ell,\text{R/B}}$. \Cref{tab:autodamp} tabulates these coefficients. Interestingly, the first level coefficients in the additive case, $\lambda_1 = \smash{\frac45}$ in 2D and $\lambda_1 = \smash{\frac67}$ in 3D, coincide with the optimal damping coefficients for Jacobi, the latter being a staple derivation in the multigrid literature \cite{mgtutorial,mgbook1,mgbook2,demmel}. As $\ell$ increases, so too does $\lambda_\ell$, suggesting that deeper cascading smoothers are less constrained by the highest frequency oscillatory error modes that are damped by prior smoothers in the cascade. %
Meanwhile, in the multiplicative case, we see that a multicoloured cascading smoother quickly transitions to overrelaxation ($\lambda > 1.1$ after the first colour sweep) and then to relatively-strong overrelaxation ($\lambda_{\ell,\text{R/B}} > 1.5$ for $\ell \geq 3$). In contrast, the optimal overrelaxation parameter for red-black Gauss-Seidel is typically\footnote{\label{footnote:rbomega}The optimal $\omega$ for red-black Gauss-Seidel is less well-defined than for the Jacobi, as it depends on various multigrid design choices, including the number of pre- and post-smoothing steps, and whether the colour ordering is reversed for pre- vs post-smoothing \cite{doi:10.1137/0917013}.} $\omega_\text{opt} \in [1, 1.2)$ (our parameter sweep yielded $\omega = 1$ in 2D and $\omega \approx 1.1$ in 3D, see \cref{fig:omega}).

\begin{figure}
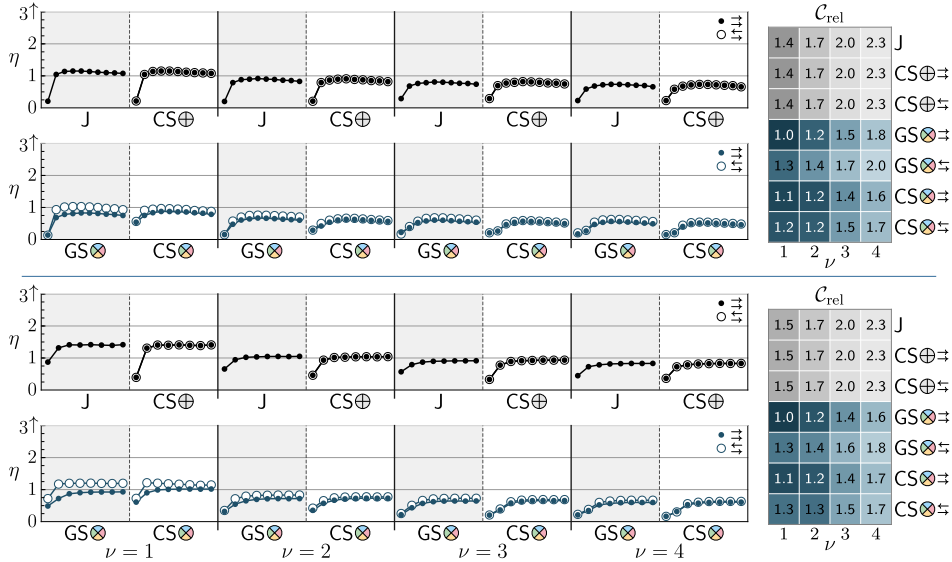
%
\centering%
\putfig{T7.Nx.P1}%
\caption{Multigrid solver performance for the 2D (top) and 3D (bottom) finite difference Poisson problem considered in \cref{sec:fd}; see \cref{sec:guide} for a visualisation guide.}%
\label{fig:T7.Nx.P1}%
\end{figure}

These insights notwithstanding, our core focus is, of course, on overall multigrid solver performance. To that end, the 2D and 3D results for this test application are shown in \cref{fig:T7.Nx.P1}.\footnote{One may note that $\mathsf{CS}\Jmode\FFmode$ yields precisely the same results as $\mathsf{CS}\Jmode\RFmode$. This atypical outcome is merely a problem-specific quirk: every $\Lambda_\ell$ in the additive cascading smoother is a scaled identity matrix, implying that the overall result is invariant to permutations as well as transpositions.} We see that a depth $\nu$ additive (resp., multicoloured multiplicative) cascading smoother yields nearly the same convergence rates as $\nu$ steps of optimally-damped Jacobi (resp., multicoloured Gauss-Seidel). In that sense, the results are not too remarkable; nevertheless, it should be emphasised that cascading smoothers do not require any parameter tuning.

\subsubsection{Finite element methods}
\label{sec:cg}

Our second test application considers a continuous Galerkin (CG) nodal finite element method (FEM) on a uniform Cartesian grid. We consider the $Q_\polydeg$ family of tensor-product polynomials of one-dimensional degree $\polydeg$, along with a Lobatto nodal basis. %
The approach we consider here is a directionally-biased per-cell non-overlapping blocking, defined by collecting the lower-left $\polydeg^d$ nodes into each cell's block. For example, for bilinear/trilinear elements, each block contains just one node; for biquadratic (resp., triquadratic) elements, a block contains 4 (resp., 8) nodes. Meanwhile, for multicoloured smoothers, we define the colouring via the associated block-adjacency graph of $A$; this leads to a 4-colour scheme in 2D and an 8-colour scheme in 3D. The multigrid interpolation operator is defined via the natural approach: fine-grid nodal values are evaluated by interpolating the input coarse-grid piecewise-polynomial function.

\begin{figure}
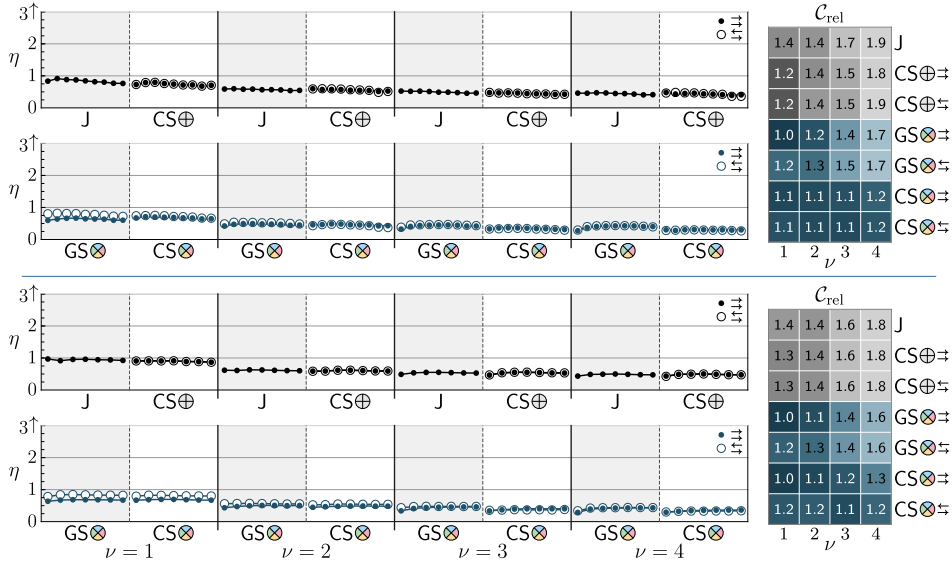
%
\centering%
\putfig{T4.Nx.P3}%
\caption{Multigrid solver performance for the 2D (top) and 3D (bottom) continuous Galerkin nodal finite element method considered in \cref{sec:cg}.}%
\label{fig:T4.Nx.P3}%
\end{figure}

We present results for the case $\polydeg = 2$ here, shown in \cref{fig:T4.Nx.P3}, and leave other cases for the \SM{}. Overall, we observe excellent multigrid performance. In the additive case, the best cascading smoother outperforms the best (optimally-damped) Jacobi method---this outcome is true almost universally, across all applications tested in this work. In the multiplicative case, note that any of the multicoloured cascading smoothers with depth $\nu \in \{1,2,3\}$ yield near-optimal cost $\mathcal{C}_\text{rel}$: this means that, as the depth increases, the convergence accelerates fast enough to offset the additional smoother cost. This could prove advantageous in more complex settings, e.g., when contending with complex geometry or stiff boundary conditions, where a depth of one is usually insufficient (as we see in following sections).

Note that the directionally biased $\polydeg^d$ blocking applied here is intended as a low polynomial order, non-overlapping baseline for comparison against classical block-Jacobi and block-Gauss-Seidel smoothers. For high polynomial orders, established $\polydeg$-independent CG FEM Poisson solvers often rely on overlapping Schwarz methods, typically using vertex-star or extended-element patches combined with a low-order coarse space or multilevel $\polydeg$-hierarchy; see, e.g., \cite{Pav93,Cas97,LF05,JK11,BF22}.
An intriguing question is whether the overlapping-block cascading smoother formulation of \cref{sec:overlapping}, built on such patches and without needing to fine-tune problem-specific damping or patch-weighting parameters, inherits the $\polydeg$-robustness of these methods. 
However, at high $\polydeg$, the corresponding Frobenius-norm minimisation problems are likely prohibitive, unless potential sketching or tensorised evaluation techniques (see \cref{sec:conclusion}) prove effective. A compelling alternative is a $\polydeg$-multigrid framework using specialised smoothers for high $\polydeg$ (such as tensorised fast-diagonalisation \cite{LF05,CGKS24}, sparse-basis Cholesky factorisation \cite{BF22}, or local $\polydeg$-multigrid patch solvers \cite{Wic25a}), combined with cascading smoother $h$-multigrid at the base level. %
We leave these directions open for future investigation.

\subsubsection{Discontinuous Galerkin methods}
\label{sec:dg}

Our next application considers a discontinuous Galerkin (DG) method, specifically a local discontinuous Galerkin (LDG) method \cite{dgmg,fluxx}, using tensor-product piecewise polynomials of one-dimensional degree $\polydeg$ and Legendre orthogonal basis. We consider uniform Cartesian grids as well as nonuniform quadtree/octree meshes, and apply the natural element-wise blocking, $(\polydeg+1)^d$ modal coefficients per block. The multigrid interpolation operator is defined via the natural approach: given a piecewise polynomial function on the coarse grid, its fine-grid interpolation is precisely the same function, but with reevaluated modal coefficients on the fine-grid elements. As in all other examples in this work, the coarse-grid operators $A_{2h}, A_{4h}, \ldots$ are equivalent to rediscretising the problem on progressively coarser meshes; for LDG, these operators can be efficiently computed in a way that obviates the need for explicit remeshing, while also mirroring the discretisation properties of the fine-grid problem, such as quadrature schemes, numerical fluxes, and penalty parameters; see \cite{dgmg,fluxx} for details.

We present two cases here, each using $\polydeg = 2$ biquadratic/triquadratic elements (and leave other cases for the \SM{}):
\begin{itemize}

\begin{figure}
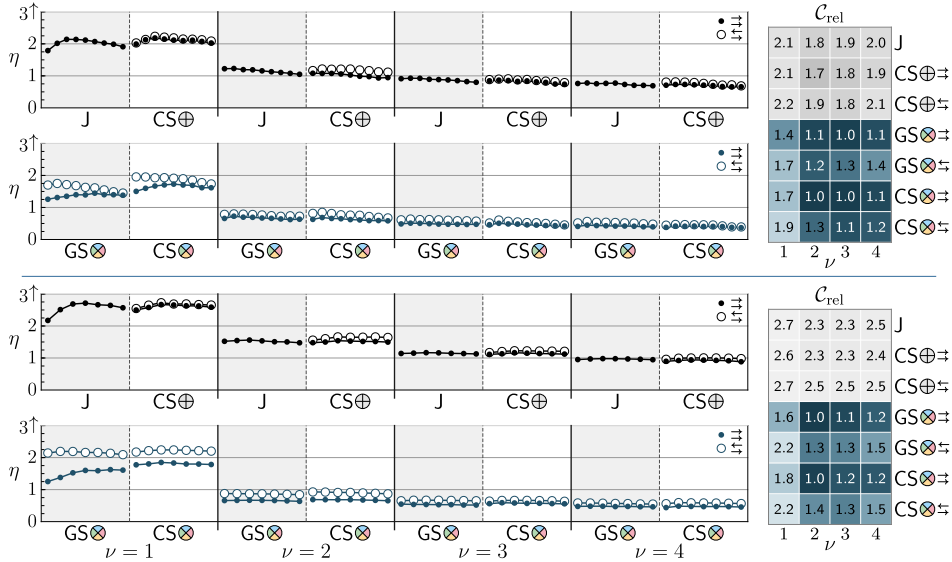
%
\centering%
\putfig{T0.Nx.P3.orbA}%
\caption{Multigrid solver performance for the 2D (top) and 3D (bottom) uniform Cartesian grid, discontinuous Galerkin problem considered in \cref{sec:dg}.}%
\label{fig:T0.Nx.P3.orbA}%
\end{figure}

    \item Uniform Cartesian grid, periodic boundary conditions --- Owing to the use of one-sided fluxes in the LDG scheme (see \cite{fluxx} for details), the block-sparse operator $A$ has a 5-point (2D) and 7-point (3D) stencil structure; accordingly, multicoloured smoothers use the canonical red-black colouring. \cref{fig:T0.Nx.P3.orbA} compiles the results. DG multigrid solvers typically require at least two pre-/post-smoothing steps to achieve the best performance, and we see that in the results: indeed, $\nu = 1$ clearly leads to suboptimal performance. Overall, the best performing cascading smoothers match, or slightly surpasses, their best performing classical counterparts. %

\begin{figure}%
\centering%
\raisebox{-0.5\height}{\includegraphics[height=1.2in]{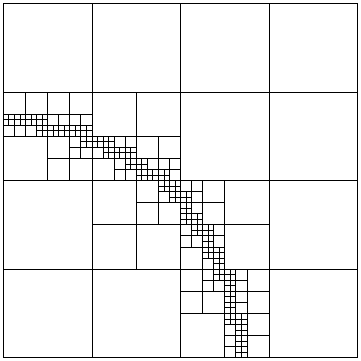}} \qquad \raisebox{-0.5\height}{\includegraphics[height=1.3in]{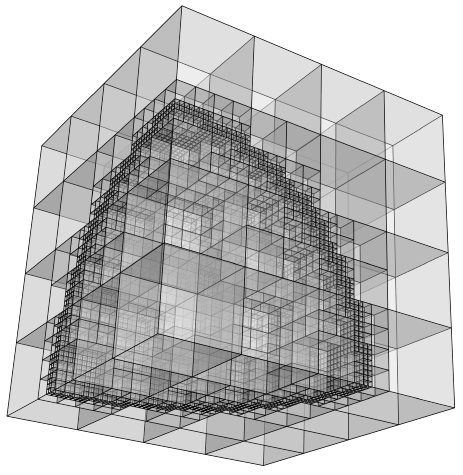}}%
\caption{Examples of the nonuniform quadtree (left) and octree (right) used in the adaptive mesh refinement problems. Each example corresponds to a base-level grid size of $n = 4$, further refined up to four levels around a circular arc or spherical shell of the same radius.}%
\label{fig:amr}%
\end{figure}

\begin{figure}
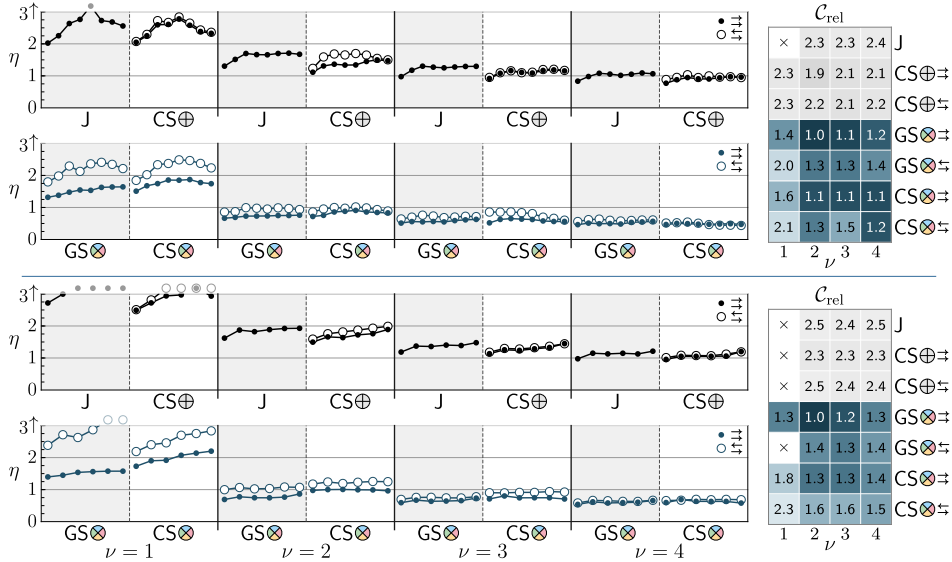
%
\centering%
\putfig{T0.Nx.P3.orbB}%
\caption{Multigrid solver performance for the 2D quadtree mesh (top) and 3D octree mesh (bottom) discontinuous Galerkin problem considered in \cref{sec:dg}.}%
\label{fig:T0.Nx.P3.orbB}%
\end{figure}

    \item Adaptively refined mesh, Neumann boundary conditions --- Illustrated in \cref{fig:amr}, the mesh has up to four levels of nested refinement around a circular arc or spherical shell. \cref{fig:T0.Nx.P3.orbB} compiles the results.\footnote{In the AMR test problems, the $n = 4, 8, 16, \ldots$ markers within each sequence correspond to the base-level grid size $n$; e.g., the last marker in 2D (resp., 3D) corresponds to $n = 512$ (resp., $128$), but the minimum element size $h$ matches that of a $8192 \times 8192$ grid (resp., $2048 \times 2048 \times 2048$). The finest-resolution problems have around 400 thousand elements in 2D and 16 million elements in 3D.} Overall, in terms of both $\eta$ and the $\nu$-weighted solver cost $\mathcal{C}_\text{rel}$, we see that a multicoloured smoother yields the fastest solvers. However, this neglects the complications of large-scale colouring: the graph colouring algorithm itself is a nontrivial component of the implementation, and if there are many smoother sweeps to perform, this can hinder massively-parallel computations. Indeed, the mesh used here is strongly nongraded (e.g., an element of width $h$ can have multiple neighbouring elements of width $h/8$), leading to a considerably complex stencil and sparsity structure, resulting in the need for $\approx 7$ colours in 2D and upwards of $\approx 40$ colours in 3D. As such, an additive-style smoother may be more effective in practical large-scale settings.
\end{itemize}

\subsection{Multiphase Elliptic Interface Problems}
\label{sec:elliptic}

Thus far, we have focused on classic Poisson problems; we next consider high-contrast multiphase elliptic interface problems, which consist of solving for a function $u: \Omega \to \R$ such that
\begin{equation} \label{eq:scalar1} \begin{aligned} -\nabla \cdot \bigl(\mu_i \nabla u \bigr) &= f && \text{in } \Omega_i, \end{aligned} \end{equation}
subject to the interfacial jump conditions
\begin{equation} \label{eq:scalar2} \left. \begin{aligned} \jump{u} &= g_{ij} \\ \jump{(\mu \nabla u) \cdot \vn} &= h_{ij} \end{aligned} \right\} \text{ on } \Gamma_{ij}, \end{equation}
and boundary conditions
\begin{equation} \label{eq:scalar3} \begin{aligned} u &= g_{\partial} && \text{on } \Gamma_D, \\ (\mu \nabla u) \cdot \vn &= h_{\partial} && \text{on } \Gamma_N, \end{aligned} \end{equation}
where $\Omega$ is a domain in $\R^d$ divided into one or more subdomains/phases $\Omega_i$, $\Gamma_{ij} := \partial \Omega_i \cap \partial \Omega_j$ is the interface between phase $i$ and $j$, and $\Gamma_D$ and $\Gamma_N$ denote the parts of $\partial \Omega$ upon which Dirichlet or Neumann boundary conditions are prescribed, respectively. Here, $\jump{\cdot}$ denotes the jump in a quantity across an interface $\Gamma_{ij}$, while $\bm n$ denotes an appropriately-oriented unit normal. Finally, $\mu_i$ is a phase-dependent ellipticity/viscosity coefficient, while $f$, $g$, and $h$ provide the data for the elliptic interface problem and are given functions defined on $\Omega$, its boundary, and internal interfaces.

We consider a high-order discontinuous Galerkin method, specifically the \textit{viscosity-upwinded} LDG methods developed in \cite{fluxx,flame}. Roughly speaking, in the limit of arbitrarily-large viscosity ratios, multiphase elliptic interface problems separate into distinct problems coupled by the interface: for the highly viscous phase, it effectively ``sees'' a Neumann-like boundary condition at the interface, whereas for the phase on the other side, it ``sees'' a Dirichlet-like boundary condition. Viscosity-upwinded methods bias the numerical fluxes of a DG method in the same way, facilitating optimal high-order accuracy.  %

\begin{figure}
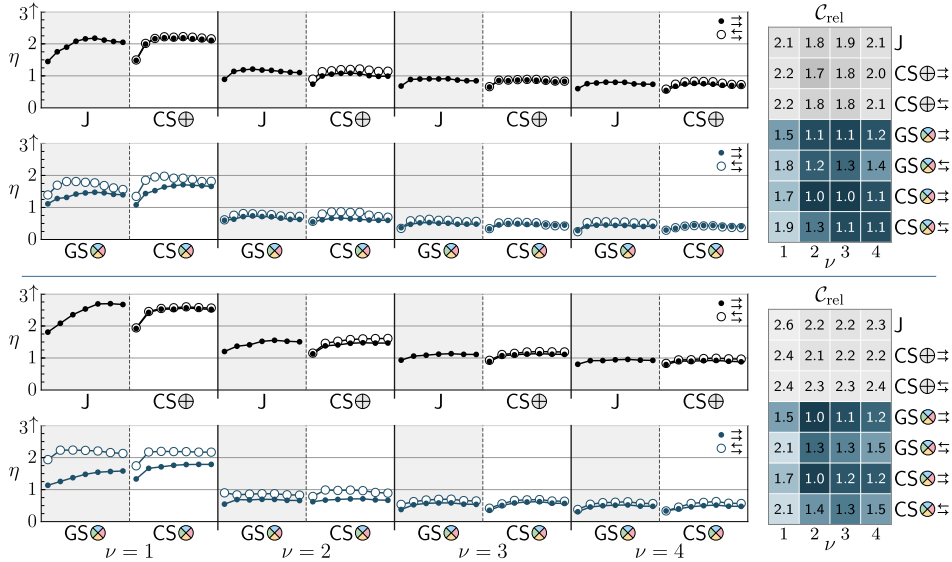
%
\centering%
\putfig{T0.Nx.P3.orbC+}%
\caption{Multigrid solver performance for the 2D (top) and 3D (bottom) high-contrast multiphase elliptic interface problem considered in \cref{sec:elliptic}.}%
\label{fig:T0.Nx.P3.orbC+}%
\end{figure}

So that we can focus on high-contrast aspects, our specific test problem considers a simplified geometry that can be easily handled by uniform Cartesian grids. Let $\Omega = (0,1)^d$ be divided into an interior rectangular phase $\Omega_0 = (\tfrac14, \tfrac34)^d$ and an exterior phase $\Omega_1 = \smash{\Omega \setminus \overline{\Omega_0}}$, each with viscosity $\mu_i$. Along with Dirichlet boundary conditions, two representative viscosity ratios are evaluated: %
$\mu_0/\mu_1 = 10^4$ and $\mu_0/\mu_1 = 10^{-4}$. Although it is important to test both extremes (as information is transferred across the domain in distinct ways), the results for both ratios are nearly identical; as such, we highlight just one configuration here, $\mu_0/\mu_1 = 10^4$ with $\polydeg = 2$, in \cref{fig:T0.Nx.P3.orbC+}, and leave other combinations to the \SM{}. Overall, these results are similar to those for a LDG method applied to a classic Poisson problem (see \cref{sec:dg}); regardless, they show that cascading smoothers are robust to high-contrast problems containing viscosity jumps many orders in magnitude; moreover, they match or slightly outperform optimally-damped Jacobi and Gauss-Seidel, but do not require any parameter tuning themselves.

\subsection{Stokes Problems}
\label{sec:stokes}

Our focus now shifts from elliptic problems to Stokes-like saddle-point problems. The general formulation consists of solving for a velocity field $\vu : \Omega \to \R^d$ and pressure field $p: \Omega \to \R$ such that
\begin{equation} \label{eq:govern1} \left. \begin{aligned} \rho_i \vu - \nabla \cdot \bigl(\mu_i (\nabla \vu + \gamma\,\nabla \vu^\trans) \bigr) + \nabla p &= \vf \\ -\nabla \cdot \mathbf \vu &= f \end{aligned} \right\} \text{ in } \Omega_i, \end{equation}
subject to the interfacial jump conditions
\begin{equation} \label{eq:govern2} \left. \begin{aligned} \jump{\vu} &= \vg_{ij} \\ \jump{\mu (\nabla \vu + \gamma\,\nabla \vu^\trans) \cdot \vn - p\, \vn} &= \vh_{ij} \end{aligned} \right\} \text{ on } \Gamma_{ij}, \end{equation}
and boundary conditions
\begin{equation} \label{eq:govern3} \begin{aligned} \vu &= \vg_{\partial} && \text{on } \Gamma_D, \\ \mu (\nabla \vu + \gamma\,\nabla \vu^\trans) \cdot \vn - p\, \vn &= \vh_{\partial} && \text{on } \Gamma_N. \end{aligned} \end{equation}
Here, $\gamma \in \{0,1\}$ is a parameter defining the form of the Stokes equations: they are said to be in \textit{standard form} when $\gamma = 0$ and \textit{stress form} when $\gamma = 1$. Most of our test applications correspond to \textit{steady-state} Stokes problems for which $\rho \equiv 0$. Conversely, if $\rho_i > 0$, then we have an \textit{unsteady} Stokes problem; the corresponding term in the momentum equations arises, e.g., through a time-stepping method applied to the incompressible Navier-Stokes equations in which $\rho_i$ is proportional to fluid density and inversely proportional to the time step.

\subsubsection{Finite difference staggered grid methods}
\label{sec:staggered}

\begin{figure}%
\centering%
\includegraphics[scale=0.9]{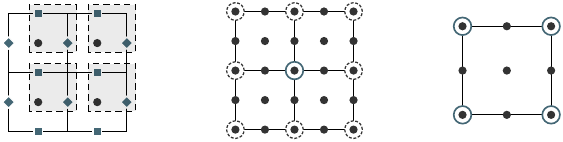}%
\newcommand{\sdiamond}{\raisebox{0.2ex}{\rotatebox[origin=c]{45}{\rule{0.8ex}{0.8ex}}}}
\newcommand{\ssquare}{\raisebox{0.2ex}{\rule{0.8ex}{0.8ex}}}
\newcommand{\opencirc}{\vcenter{\hbox{\scalebox{1.5}{$\circ$}}}}
\caption{Blocking methods for Stokes problems. %
\textbf{(left)} Finite difference staggered grid: pressure nodes are cell-centered ($\bullet$), velocity nodes are face-centered ($\sdiamond$ and $\ssquare$ for $x$- and $y$-components, respectively); dashed boxes depict the non-overlapping cell-wise blocking approach. \textbf{(middle)} Vanka patch for the Taylor-Hood element: each patch includes the central pressure node ($\opencirc$) and every velocity node ($\bullet$) from adjoining cells; neighbouring pressure nodes (dashed $\opencirc$) are not included. \textbf{(right)} Cascading smoother overlapping block approach for the Taylor-Hood element: each cell-wise block contains every node of the associated cell, $\smash{2^d}$ pressure nodes ($\opencirc$) and $\smash{3^d}$ velocity nodes ($\bullet$).}%
\label{fig:patches}%
\end{figure}

Our first Stokes test application considers a single-phase, uniform viscosity $\mu \equiv 1$, standard form, steady-state Stokes problem on a uniform Cartesian grid, discretised via the classical MAC staggered grid second-order finite difference method \cite{10.1063/1.1761178}. Pressure nodes are cell-centered, whereas velocity nodes are located at the centroid of their respective cell faces, see \cref{fig:patches}(left). One of the most effective multigrid smoothers for this problem is the Vanka smoother \cite{Vanka1986,Vanka1986b}, essentially an \textit{overlapping-block} Jacobi or Gauss-Seidel approach, where each cell-wise patch collects $2d + 1$ degrees of freedom, one for the central pressure node and $2d$ more for the surrounding velocity nodes. For example, a multiplicative-style Vanka smoother visits each cell, locally solves the Stokes system to find a coupled  update for the corresponding $2d+1$ variables, writes that update in-place, and then proceeds to the next cell; an additive-style Vanka smoother computes an update for all cells at once, and then fuses overlapping updates into a single-valued global update via weighted averaging.

In contrast, when designing smoothers for Stokes problems, the prevailing consensus is that non-overlapping block approaches lead to ineffective multigrid solvers \cite{10.1002/nla.2389}. Roughly speaking, for saddle-point systems with a constraint equation (here, velocity divergence), the blocking strategy should be such that, for every Lagrange multiplier degree of freedom within a block (pressure nodes), the block should also include all coupled primal degrees of freedom (velocity nodes); classical Vanka smoothers satisfy this criterion, but non-overlapping block methods cannot. Nevertheless, we explore the potential of both overlapping and non-overlapping block cascading smoothers:
the overlapping variant adopts the same blocking as the Vanka smoother, whereas the non-overlapping variant uses a simple cell-wise grouping shown in \cref{fig:patches}(left), collecting the central pressure node and the right-facing velocity nodes for a total of $d+1$ degrees of freedom per block.

\begin{figure}
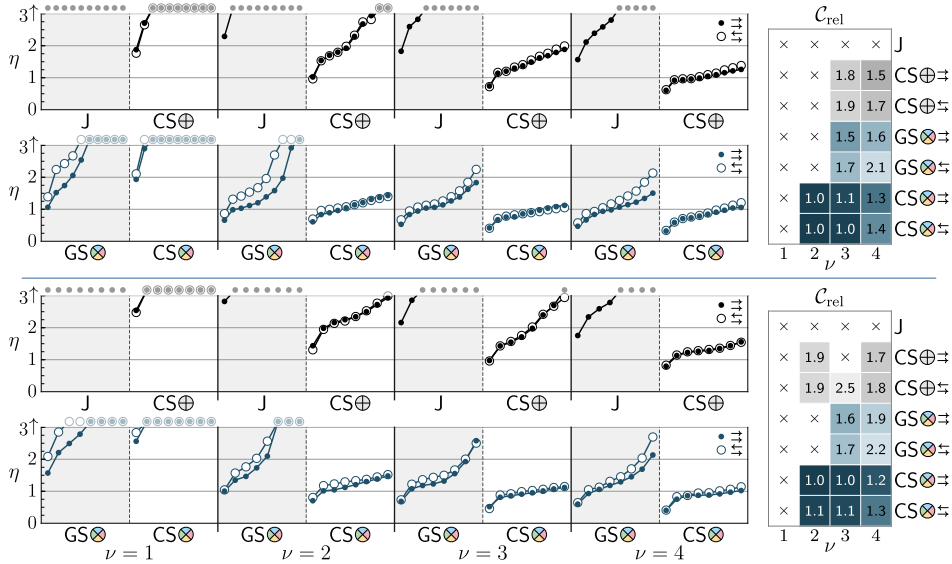
%
\centering%
\putfig{T63.Nx.P1}%
\caption{Multigrid solver performance for the 2D (top) and 3D (bottom) staggered grid finite difference Stokes problem considered in \cref{sec:staggered}.}%
\label{fig:T63.Nx.P1}%
\end{figure}

\Cref{fig:T63.Nx.P1} compiles the multigrid solver results for the staggered-grid Stokes problem with periodic boundary conditions.\footnote{The multigrid interpolation operator is standard and is defined by its action on the velocity and pressure fields: fine-grid velocity nodes are evaluated through bilinear/trilinear interpolation (appropriately accounting for the staggered grid layout); pressure nodes apply piecewise constant interpolation, i.e., the central pressure value of a coarse-grid cell is copied to all $2^d$ fine-grid subcells. Additionally, in studying \cref{fig:T63.Nx.P1}, one may observe a slow upward trend in $\eta$, even for optimally damped Vanka smoothing. This behaviour stems from a known trade-off in using the rediscretised, sparser coarse-grid operator in lieu of the denser Galerkin-projected counterpart \cite{mgbook1,mgbook2,mgbook3}.} Here, $\mathsf{V}\Jmode$ denotes the additive-style Vanka smoother, while $\mathsf{V}\GSmode\FFmode$ and $\mathsf{V}\GSmode\RFmode$ denote the red-black multiplicative variants; for all three types, the corresponding damping/overrelaxation parameter was tuned using a deep V-cycle and multigrid-preconditioned GMRES as the driver (see \cref{sec:JGS} and \cref{fig:omega}). Meanwhile, $\trimoverbar{\mathsf{CS}}$ and $\mathsf{CS}$ denote the overlapping and non-overlapping cascading smoother variants, respectively. Overall, we see that $\trimoverbar{\mathsf{CS}}$ is competitive with the otherwise analogous Vanka benchmark, yet requires no parameter fine-tuning. Intriguingly, sufficiently deep non-overlapping multiplicative cascading smoothers can yield convergence speeds $\eta$ comparable to those of their overlapping counterparts. Compared to overlapping-block approaches, non-overlapping methods tend to be simpler to implement, more cache-friendly, and more amenable to parallelisation. Combined with their smaller block sizes,\footnote{The relative solver costs $\mathcal{C}_\text{rel}$ in \cref{fig:T63.Nx.P1} do \textit{not} account for differences in block sizes. Assuming quadratic scaling in matrix-vector multiplications, the reported $\mathcal{C}_\text{rel}$ values for $\mathsf{CS}$ would be multiplied by $(3/5)^2 = 0.36$ (2D) and $(4/7)^2 \approx 0.33$ (3D); however, effective scaling is not purely quadratic at these small sizes.} multigrid solvers based on non-overlapping cascading smoothers may in practice outperform their Vanka-based counterparts, but we leave this considerably nuanced study for future work.

\subsubsection{Taylor-Hood finite element methods}
\label{sec:TH}

Our next Stokes application considers a Taylor-Hood nodal finite element method; specifically, the $\mathcal{Q}_2$--$\mathcal{Q}_1$ element consisting of biquadratic/triquadratic velocity fields and bilinear/trilinear pressure fields. We compare two families of multigrid solver, distinguished by the choice of smoother:
\begin{itemize}
    \item Vanka-type smoother --- we consider here the inclusive-patch method of \cite{10.1002/nla.2306}, illustrated in \cref{fig:patches}(b): each overlapping block/patch corresponds to a central pressure node and includes every velocity node of the adjoining cells, for a total of $d \cdot 5^d + 1$ degrees of freedom per block, i.e., 51 in 2D and 376 in 3D. %
    For all three smoother types, $\mathsf{V}\Jmode$, $\mathsf{V}\GSmode\FFmode$, and $\mathsf{V}\GSmode\RFmode$, optimal damping/overrelaxation was determined using a deep V-cycle and multigrid-preconditioned GMRES as the driver (see \cref{sec:JGS} and \cref{fig:omega}).
    \item Overlapping-block cascading smoother --- illustrated in \cref{fig:patches}(c), we consider a cell-wise blocking approach, i.e., each block corresponds to the set of nodal values on each cell, for a total of $2^d + d\cdot3^d$ degrees of freedom per block, i.e., 22 in 2D and 89 in 3D. Note that this blocking strategy runs counter to the prevailing wisdom on effective blocking methods for Stokes (see \cref{sec:staggered}), because a block contains pressure nodes without including all of their coupled velocity nodes (specifically, those from neighbouring cells).
\end{itemize}

\begin{figure}
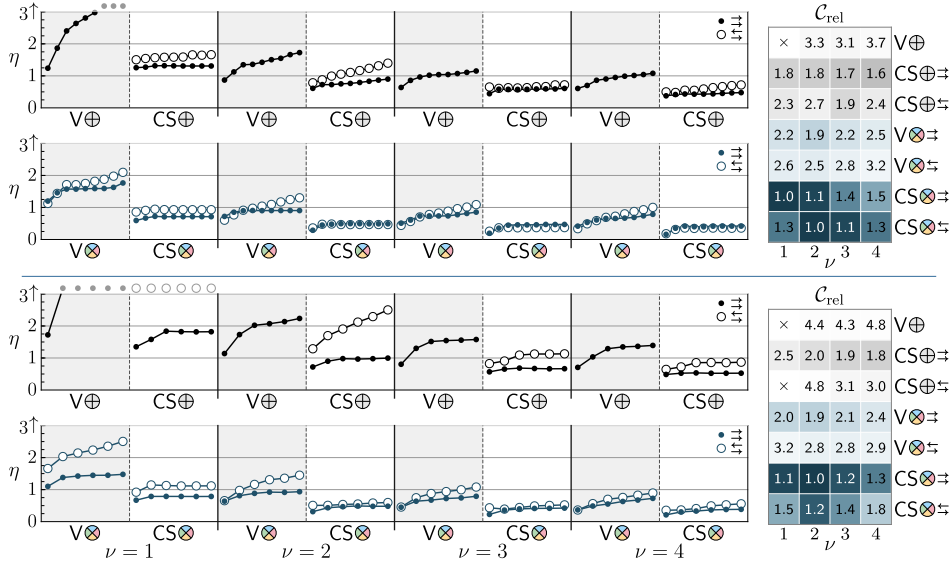
%
\centering%
\putfig{T5a.Nx.P3}%
\caption{Multigrid solver performance for the 2D (top) and 3D (bottom) Taylor-Hood finite element Stokes problem considered in \cref{sec:TH}.}%
\label{fig:T5a.Nx.P3}%
\end{figure}

\noindent With this setup, \cref{fig:T5a.Nx.P3} compiles the corresponding multigrid solver results. Overall, we see that the cascading smoother approach leads to remarkably effective multigrid solvers:\footnote{In 3D, $\mathsf{CS}\GSmode\RFmode$ performs poorly, with $\eta$ deteriorating under grid refinement; however, rapid convergence can be reclaimed by using an alternative prescaling norm---see the \SM{}.} at their respective optimal configurations of $\nu$ and $\RFmode$/$\FFmode$, $\mathsf{CS}\Jmode$ (resp., $\mathsf{CS}\GSmode$) converges about twice as fast as $\mathsf{V}\Jmode$ (resp., $\mathsf{V}\GSmode$).\footnote{Various weighting strategies can be used for Vanka smoothers \cite{10.1002/nla.2306}. Our implementation applies node-dependent weighted averaging, multiplied by a global damping factor $\omega$; e.g., velocity nodes shared by $m$ cells apply a weight of $\omega/m$, whereas pressure nodes have weighting $\omega$ (there is only one pressure node per block). This approach is applied for both additive-style and multiplicative-style Vanka smoothers as it led to the fastest convergence in both formulations. Although we did not comprehensively explore other weighting strategies, the convergence rates shown here are representative of tuned Vanka solvers for Taylor-Hood elements (see \cite{10.1002/nla.2306}).} Further, it should be emphasised that the relative solver costs $\mathcal{C}_\text{rel}$ shown in \cref{fig:T5a.Nx.P3} do \textit{not} account for the difference in block sizes. Relative to the blocking of the cascading smoother approach, the Vanka method is a factor $\approx 2.3$ bigger in 2D and $\approx 4.2$ in 3D; since matrix-vector multiplications scale approximately quadratically in block size, it is reasonable to multiply the reported Vanka $\mathcal{C}_\text{rel}$ factors by $\approx 5$ and $\approx 18$, respectively. This reinforces our expectation that a cascading smoother multigrid solver offers significant potential for Taylor-Hood finite element methods, and perhaps other CG-based Stokes frameworks as well.

\subsubsection{Discontinuous Galerkin methods}
\label{sec:DGstokes}

Our last set of test applications considers a mixed-degree discontinuous Galerkin method, specifically a local discontinuous Galerkin (LDG) framework, initially developed by Cockburn \textit{et al} \cite{CockburnKanschatSchotzauSchwab2002} and extended to multiphase variable-viscosity problems in \cite{thunderball}. In particular, we consider a tensor-product Legendre polynomial orthogonal basis of one-dimensional degree $\polydeg$ (resp., $\polydeg - 1$) for the velocity field (resp., pressure field). The multigrid interpolation operator is defined via the natural approach: given a piecewise polynomial velocity/pressure field on the coarse grid, its fine-grid interpolation is precisely the same function, but with reevaluated modal coefficients on the fine-grid elements. Once again, the coarse-grid Stokes problem is equivalent to rediscretisation; notably, the corresponding operators $A_{2h}, A_{4h}, \ldots$ can be efficiently computed in a way that obviates the need for explicit remeshing, while also mirroring the fine-grid discretisation with regard to properties such as quadrature schemes, numerical fluxes, and penalty parameters; see \cite{flame} for details on the corresponding operator-coarsening methods. %

\begin{figure}%
\centering%
\putfig{T2.Nx.P3.orbAbcP.CS}%
\caption{Corresponding to a $\polydeg = 2$ mixed-degree LDG method, multigrid solver performance for a 2D (left) and 3D (right) single-phase, uniform viscosity, steady-state, stress-form Stokes problem on a uniform Cartesian grid with periodic boundary conditions.}%
\label{fig:T2.Nx.P3.orbAbcP}%
\end{figure}

In the setting of a DG method, it is especially natural to consider an element-wise blocking strategy, whereby each block corresponds to the collective set of velocity and pressure degrees of freedom on each mesh element; this non-overlapping approach is used here. Unlike prior test problems, which focused on canonical standard-form Stokes systems, the DG formulation accommodates a wider range of boundary conditions and enables us to explore more challenging stress-form problems.

Our first example considers a single-phase, uniform viscosity, steady-state, stress-form Stokes problem on uniform Cartesian grids with periodic boundary conditions; results are shown in \cref{fig:T2.Nx.P3.orbAbcP}. We omit comparisons to standard element-wise block Jacobi and block Gauss-Seidel smoothers, as they prove to be ineffective for mixed-degree LDG Stokes problems, even with the best possible damping \cite{thunderball}.\footnote{Comparisons with classical smoothers are retained in the \SM{} for completeness.} Note that depth $\nu = 1$ cascading smoothers are also ineffective; roughly speaking, they may need at least two levels to adapt to the saddle-point structure of the discretised Stokes problem. From \cref{fig:T2.Nx.P3.orbAbcP} we see that a depth 2--3 (resp., 3--4) cascading smoother yields the best-performing multiplicative (resp., additive) variant; similar outcomes hold for other problem configurations including polynomial degrees, boundary conditions, and adaptively refined meshes, as shown in the \SM{}. For comparison, in \cite{thunderball} the author developed a specialised smoother for mixed-degree LDG Stokes, analogous to element-wise block Gauss-Seidel, except that the diagonal block inverses are replaced with a three-parameter approximation balancing velocity and pressure smoothing. Once tuned, the corresponding multigrid solver (calibrated solely to the equivalent of $\nu = 3$) yields $\eta \approx 0.8$ in 2D and $\eta \approx 0.9$ in 3D \cite[Fig.\,5]{thunderball}. This is about 2--10\% faster than the $\mathsf{CS}\GSmode\RFmode_{\nu=3}$ results shown in \cref{fig:T2.Nx.P3.orbAbcP}, both in terms of $\eta$ and per-iteration smoother runtime. However, the cited method requires an expensive, LFA-based three-parameter optimisation that must be repeated for each distinct problem configuration (i.e., spatial dimension, polynomial degree, standard- or stress-form, and steady or unsteady regime); furthermore, the approach in \cite{thunderball} only considered multicoloured smoothers. In contrast, cascading smoothers require significantly less tuning and naturally extend to Jacobi-style additive variants.

\begin{figure}
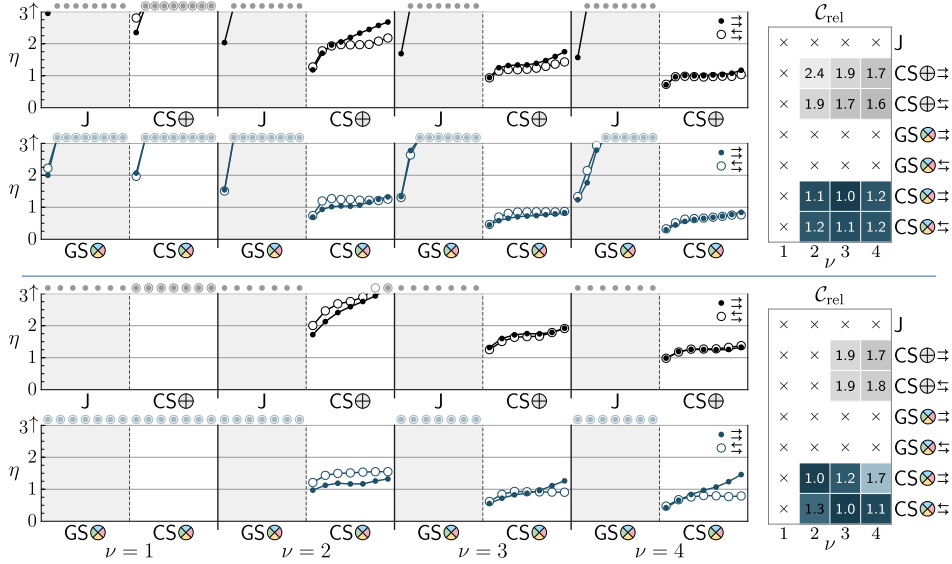
%
\centering%
\putfig{T2.Nx.P3.orbC-.CS}%
\caption{Corresponding to a $\polydeg = 2$ mixed-degree LDG method, multigrid solver performance for a 2D (left) and 3D (right) high-contrast multiphase, steady-state, stress-form Stokes problem on a uniform Cartesian grid with Dirichlet boundary conditions.}%
\label{fig:T2.Nx.P3.orbC-}%
\end{figure}

Next we consider a high-contrast multiphase problem. %
Let $\Omega = (0,1)^d$ be divided into an interior rectangular phase $\Omega_0 = (\tfrac14, \tfrac34)^d$ and an exterior phase $\Omega_1 = \Omega \setminus \overline{\Omega_0}$, each with viscosity $\mu_i$. We consider a viscosity ratio $\mu_0/\mu_1 = 10^{-4}$, with other ratios yielding similar results (see the \SM{}). Applying Dirichlet boundary conditions (i.e., $\Gamma_N = \varnothing$ and $\Gamma_D = \partial \Omega$ in \cref{eq:govern3}), \cref{fig:T2.Nx.P3.orbC-} compiles the multigrid solver results. Once again, cascading smoothers yield compelling results, with depth 2 or 3 smoothers providing the fastest wall-clock time. One aspect visible in the 3D results is that $\mathsf{CS}\GSmode\FFmode$ degrades in performance as the grid is refined, but $\mathsf{CS}\GSmode\RFmode$ yields grid-independent convergence rates. This example shows that, although cascading smoothers do not require parameter fine-tuning, in practice we may still need to test which forward/reverse application yields the best results.

\begin{figure}
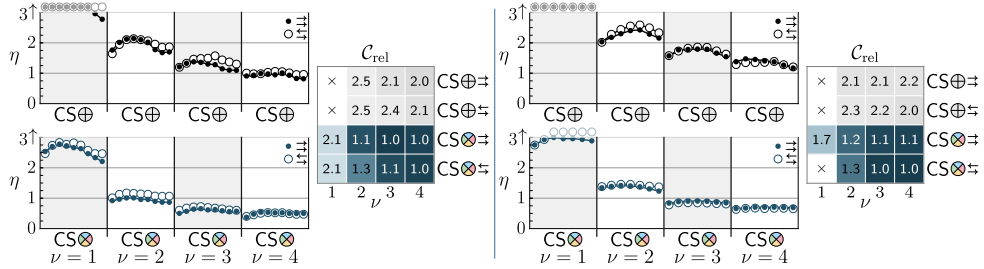
%
\centering%
\putfig{T2.Nx.P3.orbD2.CS}%
\caption{Corresponding to a $\polydeg = 2$ mixed-degree LDG method, multigrid solver performance for a 2D (left) and 3D (right) single-phase, unsteady stress-form Stokes problem on a uniform Cartesian grid with periodic boundary conditions.}%
\label{fig:T2.Nx.P3.orbD2.CS}%
\end{figure}

Our final application serves to highlight the nuances of unsteady Stokes problems. The introduction of a density-weighted temporal-derivative term in the Stokes momentum equations can alter the characteristics of the saddle-point problem and thus the multigrid solver. In one extreme, when $\rho$ is sufficiently small, we essentially have a perturbed steady-state Stokes problem with a small $\rho$-weighted identity shift added to the viscous operator. In the other extreme, the viscous operator is essentially negligible and \cref{eq:govern1} approximately reduces to a flux-form Poisson problem for pressure. The role of pressure changes in these two extremes, from being a Lagrange multiplier enforcing the divergence constraint (in the former) to being the primary solution variable in an elliptic interface problem (in the latter). Crucially, the relative strengths of these operators can also change across the multigrid hierarchy itself---on highly refined meshes, viscous effects may dominate, whereas on coarse grids, the density-weighted term may dominate. As such, a multigrid method may have to adapt its underlying parameters depending on the regime; in prior work, the author tackled this issue through a harmonic-weighted blending scheme that interpolated between two extremal sets of optimally-tuned solver parameters \cite{flame,thunderball}. On the other hand, a cascading smoother-based multigrid solver, by design, has no parameters to tune; consequently, we wonder whether it can handle a broad range of effective Reynolds numbers. To highlight one such test, we consider a moderately-weak viscous problem with $\mu = 10^{-4}$ and $\rho = 1/(0.1 h)$, where $h$ is the cell size of the primary mesh; in effect, $\rho$ represents a unit-density temporal-derivative term arising from a time-stepping method with unit order CFL. \Cref{fig:T2.Nx.P3.orbD2.CS} contains the corresponding multigrid solver results. We observe reliable convergence with solver speeds $\eta$ typical of LDG-based Stokes problems; the same outcome holds for other parameter regimes as well as for unsteady high-contrast multiphase problems, as shown in the \SM{}. As a final remark, cascading smoothers were also tested and found to be effective in an equal-degree LDG Stokes framework (in which velocity and pressure fields are of equal polynomial degree, alongside pressure-penalty stabilisation).

\section{Concluding Remarks}
\label{sec:conclusion}

Cascading smoothers offer a compelling and versatile addition to the multigrid solver landscape, proving highly effective across a broad range of problems. By design, they are parameter-free except for a few discrete choices relating to their forward or reverse application in the pre- or post-smoothing step, and in the case of Stokes problems, whether to prescale the divergence constraint using the 1-norm or the max-norm (see \cref{sec:prescale}). Of the two ordering approaches, i.e., $\RFmode$ and $\FFmode$ (see \cref{sec:vcycle}), we saw that in many cases the difference in convergence rate is marginal; DG-based problems showed a greater dependence, although it is straightforward to choose the optimal approach.

Classical smoothers, such as Jacobi, Gauss-Seidel, Vanka, and Chebyshev polynomial smoothers, are often tuned in an idealised setting, e.g., uniform grids with periodic boundary conditions or via two-grid local Fourier analysis \cite{mgbook1,mgbook2,doi:10.1137/19M1308669}. As such, their performance can be impacted by boundary conditions, stiff penalty parameters, highly unstructured grids, complex geometry, etc. For example, polynomial smoothers are vulnerable to isolated extremal eigenvalues outside the bulk spectrum. In contrast, cascading smoothers are constructed locally, optimising for the best sequence of block-diagonal smoothers that minimises the Frobenius norm of the cascading error propagator, taking into account the local geometry, local ellipticity coefficients, nearby boundaries and interfaces, etc. Our results show they are robust to different kinds of boundary conditions, strongly nongraded adaptively refined meshes, and coefficients that vary over multiple orders magnitude.

As discussed in \cref{sec:impl}, progressively-deeper cascading smoothers are more expensive to build owing to residual fill-in. For many applications and/or shallow cascading smoothers this is not an issue, e.g., whenever the smoother can be precomputed offline. Several approaches could be taken to accelerate their construction. A simple method is to limit fill-in by dropping all but the largest $\alpha \in \mathbb{N}$ blocks of the smoother residual (e.g., as measured by the Frobenius norm); in many of the applications considered here, $\alpha \approx 5$ in 2D and $\alpha \approx 7$ in 3D yields practically-similar multigrid solver speeds as $\alpha = \infty$. Another technique might take advantage of stencil repetition on structured meshes: e.g., a uniform-coefficient Poisson/Stokes problem on a uniform Cartesian grid yields the same smoother construction almost everywhere across the domain; those calculations could be cached and reused. Beyond these approaches, a promising reformulation is to make use of sketching methods. For an additive cascading smoother $(\Lambda_\ell)_{\ell=1}^\nu$, recall that $\Lambda_\ell := \argmin_{\Lambda \in \D} \|(I - \Lambda A)E_{\ell - 1}\|_F$, where $E_\ell := (I - \Lambda_\ell A)E_{\ell - 1}$ is the residual at each level and $E_0 := I$. A sketched approach would replace the square minimisation problem %
with a tall-and-narrow rectangular approximation $\min_\Lambda \|(I -\Lambda A) E_{\ell-1} S\|_F$, where $S$ is sketching matrix (e.g., Gaussian or Rademacher). Unwinding the recursion, this is equivalent to the following: imagine $S$ to be a list of randomly generated error modes; then, find the best possible block-diagonal smoother that damps all of them at once (in the sense of the total squared $L^2$ norm residual); replace $S$ with the smoothed result, $S \leftarrow (I - \Lambda_1 A)S$, and repeat the process to find $\Lambda_2$, etc. Note that the sketched problem remains trivially parallelisable over the set of diagonal blocks of $\Lambda_\ell$. The formulation naturally extends to multicoloured cascading smoothers and may also lend itself to matrix-free implementations (e.g., those where the system operator $A$ is implemented as a black-box). Key questions include: what sketching size is needed for near-optimal results; whether the sketched problems and resulting smoothers are robust to limited input randomness; %
and whether the approach is competitive with deterministic truncation/dropoff methods or if a hybrid approach could be used (e.g., direct solve on the first level, sketch on remaining levels). We leave this exploration open for a future study.

Looking ahead, we anticipate that cascading smoothers will find application across a rich array of problems, well beyond the scope of this initial study. For example, different multigrid designs could be explored, e.g., cascading smoothers in combination with F-cycles or W-cycles as well as pre- and post-smoothers of unequal depths. For very high-order problems, as discussed at the end of \cref{sec:cg}, cascading-smoother $h$-multigrid could be used as a general-purpose base solver within a specialised $\polydeg$-robust framework; alternatively, incorporating tensor-product structure or sketching could mitigate the $\mathcal{O}(\polydeg^{3d})$ explicit construction cost at high order. Besides the quadtree/octree examples of this work, cascading smoothers should naturally extend to problems with unstructured simplicial meshes as well as complex geometry; e.g., preliminary work indicates they are highly effective on implicitly-defined cut cell meshes containing abrupt and extreme variations in element shape and size. Additionally, an intriguing open question is whether the cascading smoother concept can be directly applied or otherwise adapted to problems with strong directional biasing, such as convection-dominated or anisotropic problems as well as other kinds of heterogeneity, such as in composite materials or porous media. Finally, although we focused here on geometric multigrid methods, cascading smoothers are fundamentally agnostic to the specific grid-transfer operators; coupling them to AMG methods is another compelling avenue, leveraging locally adaptive smoothing alongside AMG's strengths in adaptive coarse-space generation.

\section*{Acknowledgements}

This work was supported by a U.S. Department of Energy, Office of Science Early Career Research Program award and by the U.S. Department of Energy, Office of Science, Office of Advanced Scientific Computing Research's Applied Mathematics Competitive Portfolios program, each under Contract No.~DE-AC02-05CH11231. Some computations used resources of the National Energy Research Scientific Computing Center (NERSC), a U.S. Department of Energy, Office of Science User Facility operated at Lawrence Berkeley National Laboratory.

\bibliographystyle{siamplain}
\bibliography{references}

\clearpage

\fi

\ifIncludeSM

\begin{center}
{\bfseries CASCADING SMOOTHERS FOR MULTIGRID \\[1ex] \MakeUppercase{\SMx}\par}%
\vskip2mm%
{\small\textsc{Robert I. Saye \orcidlink{0000-0001-7256-6941}}\par}%
\vskip1mm%
{\sffamily\footnotesize Lawrence Berkeley National Laboratory, Berkeley, California, USA $\mid$ \texttt{rsaye@lbl.gov} $\mid$ \paperdate\par}%
\vskip5mm%
\end{center}

\renewcommand{\thesection}{\scalebox{0.7}{\SM{}}\,\arabic{section}}
\renewcommand{\thefigure}{\scalebox{0.7}{\SM{}}\,\arabic{figure}}
\renewcommand{\thetable}{\scalebox{0.7}{\SM{}}\,\arabic{table}}

\section{Multicoloured Cascading Smoothers}

\subsection{Nonoverlapping Blocks} 

As noted in the main article, a multiplicative cascading smoother can be built for any kind of colouring scheme. However, if the scheme corresponds to a graph colouring of the block-sparsity structure of $A$, then we can simplify various implementation aspects by overwriting state variables in-place, analogous to classical multicoloured Gauss-Seidel methods. Partitioning $A$ into elementwise blocks $A_{ij}$, suppose that we have coloured the grid such that elements $i$ and $j$ have different colours whenever $A_{ij}$ is nonzero and $i \neq j$.\footnote{Minimising the number of distinct colours helps to maximise parallel efficiency and minimise fill-in of $E$. However, it is not necessary to perfectly colour the system; a near-optimal colouring remains highly effective.} One can then replace the level-and-colour-dependent residuals $E_{\ell,k}$ with a single block-sparse matrix $E$ that is updated during construction, as shown in \cref{algo:SM-mmcs}. Note that multiple block rows of $E$ can be updated concurrently (line \ref{algo:SM-update}); this does not affect any other computation in the same sweep, e.g., the calculation on line \ref{algo:SM-soln} only requires data from prior sweeps.

\begin{figure}[!h]%
\centering\vspace{-1em}%
\begin{minipage}{11.2cm}%
\begin{algorithm}[H]
	\caption{Building a depth $\nu$ multicoloured cascading smoother $(\Lambda_{\ell,k})_{\ell=1,k=1}^{\nu,\chi}$.}
	\begin{algorithmic}[1]
        \State Compute an appropriate diagonal prescaling matrix $W$.
        \State Prescale $A$ by replacing $A \leftarrow W\!A W$.
        \State Initialise $E := \Id$.
        \For{level $\ell = 1, \ldots, \nu$}
            \For{each colour $k = 1, \ldots, \chi$}
                \For{every element $i$ of colour $k$ (possibly in \textbf{parallel})} \label{line:coloursweep}
                    \State Define $\Lambda_{\ell,k,ii}$ as the least squares solution of
                    \Statex \qquad $\Lambda_{\ell,k,ii} (A E)_{i,:} = E_{i,:}.$ \label{algo:SM-soln}
                    \State Update the $i$th block row of the residual,
                    \Statex \qquad $E_{i,:} \leftarrow E_{i,:} - \Lambda_{\ell,k,ii} (A E)_{i,:}.$ \label{algo:SM-update}
                    \State Optionally, truncate the residual, e.g., by removing all but the largest $\alpha$ blocks of $E_{i,:}$.
                    \State Transform to unscaled system,
                    \Statex \qquad $\Lambda_{\ell,k,ii} \leftarrow W_{ii} \Lambda_{\ell,k,ii} W_{ii}.$
                \EndFor
            \EndFor
        \EndFor
        \State Restore $A$ to its unscaled form, $A \leftarrow W^{-1}\!A W^{-1}$.
	\end{algorithmic}%
	\label{algo:SM-mmcs}%
\end{algorithm}%
\end{minipage}%
\end{figure}

\subsection{Overlapping Blocks}

The implementation of an overlapping-block multicoloured cascading smoother is essentially identical to \cref{algo:SM-mmcs}, except that, rather than using a single fixed block partition, $A$ is implicitly repartitioned for each colour immediately prior to the innermost loop. Within the loop, the residual $E$ is updated for all rows participating in the sweep of colour $k$, while all other parts of $E$ remain unchanged. Because of overlap, an individual row of $E$ can be updated multiple times per outer level $\ell$. Post-construction, application of the smoother is algorithmically analogous to an overlapping-block Gauss-Seidel method.

Algebraically, the formulation of an overlapping-block multiplicative cascading smoother---whether properly or improperly coloured---is identical to that given in the main article, except for a modified $\Lambda$ solution space. Explicitly, suppose we have $\chi$ different colours. For each colour $k \in \mathbb{N} \cap [1,\chi]$, let $\mathbb{D}_k$ denote the effective block-diagonal matrix solution space: if $\Lambda \in \mathbb{D}_k$ and $\Lambda_{ij} \neq 0$, then $i$ and $j$ must belong to the same block and that block must be of colour $k$. Assuming $A$ has been appropriately prescaled, then $\Lambda_{\ell,k} := \argmin_{\Lambda \in \D_k} \| (\Id - \Lambda A) E_{\ell,k-1} \|_F,$ where $E_{\ell,k} := (\Id - \Lambda_{\ell,k} A) E_{\ell,k-1}$ and $E_{\ell,0} := E_{\ell-1,\chi}$ if $\ell > 1$ else $E_{1,0} := \Id$. Given an input approximation $x_0$ and a right-hand side $b$, a single (forward) application is defined by the output $x_{\nu,\chi}$ where $x_{\ell,k} := x_{\ell,k-1} - \Lambda_{\ell,k} (A x_{\ell,k-1} - b)$ and $x_{\ell,0} := x_{\ell - 1, \chi}$ if $\ell > 1$ else $x_{1,0} := x_0$.

\section{Prescaling Stokes Problems}
\label{sec:SM-prescale}

In general, Stokes problems take the form
\begin{equation} \label{eq:SM-general} Ax = b \quad \Leftrightarrow \quad \begin{pmatrix} \mathcal{M}_\rho + \mathcal{A}_\mu & \mathcal{G} \\ \mathcal{D} & \mathcal{P} \end{pmatrix} \begin{pmatrix} \vu_h \\ p_h \end{pmatrix} = \begin{pmatrix} b_{\vu} \\ b_{p} \end{pmatrix}\!, \end{equation}
where $\mathcal{M}_\rho$ is a density-weighted mass matrix, $\mathcal{A}_\mu$ is a viscosity-weighted Laplacian-like viscous operator, $\mathcal{G}$ is a pressure gradient operator, and $\mathcal{D} = \mathcal{G}^\trans$ is a negative-divergence operator. These operators locally scale as $\mathcal{M}_\rho \sim \rho h^d$, $\mathcal{A} \sim \mu h^{d-2}$, $\mathcal{G} \sim h^{d-1}$, and $\mathcal{D} \sim h^{d-1}$.\footnote{If $\mathcal{P}$ is nonzero, its scaling is more complex and relates to the Schur complement; e.g., for an equal-degree DG method, depending on various discretisation choices, $\mathcal{P} \sim h^d (\mu + \rho h^2)^{-1}$.}
As such, their relative strength can change across the multigrid hierarchy. For the purposes of building a cascading smoother, the system is temporarily prescaled via a simple symmetric diagonal scaling, replacing $A$ with $W\!AW$, where $W := \diag(w_\vu,w_p)$. %
The goal is to design a scaling approach that is effective across the entire parameter regime and across the multigrid hierarchy, from steady-state Stokes problems (corresponding to $\mathcal{M}_\rho = 0$ and $\mathcal{A}_\mu \succeq 0$), to vanishing-viscosity unsteady Stokes problems ($\mathcal{M}_\rho \succ 0$ and $\mathcal{A}_\mu = 0$). %
The approach used here scales the individual blocks of \cref{eq:SM-general} by first scaling the top-left symmetric positive (semi)definite block $\mathcal A := \mathcal{M}_\rho + \mathcal{A}_\mu$, similar to an elliptic problem, after which the rows of the bottom-left block (equiv., columns of the top-right block) are normalised. Specifically:
\begin{enumerate}
    \item The top-left block is scaled so that it has unit diagonals weighted by the relative $L^1$ norm of its rows, i.e., so that the $i$th diagonal of $\diag(w_\vu) \mathcal{A} \diag(w_\vu)$ is equal to $\| \mathcal{A}_{i,:}\|_1 / \mathcal{A}_{ii}$; explicitly, $w_{\vu,i} := \mathcal{A}_{ii}^{-1} \sqrt{\| \mathcal{A}_{i,:}\|_1}$.
    \item Next, $w_p$ is chosen to unitise the $L^q$ norm of the rows of the bottom-left block, i.e., so that each row of $\diag(w_p) \mathcal{D} \diag(w_\vu)$ has $L^q$ norm equal to $1$.
\end{enumerate}
This leaves one parameter to choose, $q$. We found that it is unnecessary to perform a finely-tuned parameter sweep; instead, one can simply choose from $q \in \{1,2,\infty\}$, say. Besides being relatively simple to implement, the $L^1$ weighting strategy in part (i) has two notable features: first, the method is agnostic to the discretisation framework; second, the vanishing-density regime and vanishing-viscosity regime have nearly the same optimal $q$ (thus, one can prioritise the vanishing-density regime, say, and use that $q$ independent of the effective Reynolds number). Roughly speaking:
\begin{itemize}
    \item In the inviscid case ($\mu \equiv 0$), scaling $\mathcal{A}$ has the same effect as normalising the density. For example, for finite difference methods and orthogonal-basis DG methods, $\mathcal{A}$ is a density-scaled identity matrix, the relative $L^1$-norm weighting strategy in (i) has no effect, and the prescaled form of $\mathcal{A}$ is simply the identity matrix. The role of the $q$-norm scaling is to remove the scaling in $h$ and balance the strength of the $\mathcal{D}$ and $\mathcal{G}$ blocks relative to the top-left identity block, the latter acting as a conduit between primary solution and flux variables.
    \item In the steady-state case ($\rho \equiv 0$), the top-left block is a Laplacian-like operator coupling significantly more system variables. Separating out the scaling dependence on $\mu$ and $h$, the strength of this coupling should be approximately preserved, so as to balance the smoothing of the velocity and pressure variables. The relative $L^1$ norm weighting strategy of (i) fulfils this function, and is closely related to estimating the maximal eigenvalues of the highest-frequency eigenfunctions.
\end{itemize}
Coupled with cascading smoothers, this scaling approach is remarkably effective across a multitude of Stokes discretisation frameworks, including staggered grid finite difference schemes, Taylor-Hood finite element methods, mixed-degree LDG methods, and equal-degree LDG methods. In the finite difference case, $q = 1$ led to the most effective cascading smoother-based multigrid solvers, whereas in all other cases, $q = \infty$ is best.\footnote{It is possible to fine-tune this approach. For example, we could scale the rows of $\mathcal{D}$ to have a max-norm of $\theta \in \R_+$, say, and then meticulously tune $\theta$ for maximal multigrid convergence. Brief numerical experiments indicate potential multigrid solver speedups, albeit minor (e.g., 10\% faster). The optimal parameter would depend on a plethora of application-specific factors, including the discretisation method, polynomial degree, spatial dimension, stress-form or standard-form, steady-state or vanishing-viscosity, etc. Clearly, tuning $\theta$ for each application comes at great expense and sacrifices the elegance of the near-universal approach.} We hypothesise this is related to the sparsity of the discretisation stencils: low-order finite difference stencils are relatively flat, whereas higher-order stencils have sharp peaks; the $1$-norm (resp., max-norm) best measures their ``energy'', at least for the purposes of balancing the least squares construction underlying cascading smoothers.

\section{V-Cycle Symmetry}
\label{sec:Vcyclesymmetry}

In this section we analyse the symmetry of a cascading smoother V-cycle, using the reverse application as pre-smoother and forward application as post-smoother. While the inductive proof of V-cycle symmetry under adjoint pre- and post-smoothing is a classical result (see, e.g., Bank \& Douglas, SIAM J.~Numer.~Anal., 1985; Hackbusch, Multi-Grid Methods and Applications, Springer, 1985; and Bramble \& Pasciak, Math.~Comp., 1987), a self-contained presentation is given here to explicitly account for the forward and reverse applications inherent to cascading smoothers.

Let $b$ be the given right-hand side. Applied to an arbitrary input $x_0$, one can show (e.g, via induction) the forward application of a cascading smoother $(\Lambda_\ell)_{\ell=1}^\nu$ yields the output
\[ \mathsf{forward} (\Lambda) : x_0 \mapsto (I - \Lambda_\nu A) \cdots (I - \Lambda_1 A) x_0 + A^{-1} \Bigl[ I - (I - A \Lambda_\nu) \cdots (I - A \Lambda_1) \Bigr] b. \]
(If $A^{-1}$ does not exist, the bracketed expression is to be interpreted symbolically in the sense that every term of the expanded form of $I - (I - A \Lambda_\nu) \cdots (I - A \Lambda_1)$ has the left-most $A$ struck off.) Similarly, the reverse application yields the output
\[ \mathsf{reverse} (\Lambda) : x_0 \mapsto (I - \Lambda_1^\trans A) \cdots (I - \Lambda_\nu^\trans A) x_0 + A^{-1} \Bigl[ I - (I - A \Lambda_1^\trans) \cdots (I - A \Lambda_\nu^\trans) \Bigr] b. \]
With this setup, the V-cycle proceeds as follows. Starting with an input approximation $x_0 = 0$, apply the pre-smoother to obtain the updated solution
\[ x_\text{pre} := A^{-1} \Bigl[ I - (I - A \Lambda_1^\trans) \cdots (I - A \Lambda_\nu^\trans) \Bigr] b. \]
Next, calculate the residual
\[ r := b - A x_\text{pre} = (I - A \Lambda_1^\trans) \cdots (I - A \Lambda_\nu^\trans) b. \]
Restrict the residual, apply the coarse-grid V-cycle $V_{2h}$, and interpolate the correction:
\[ x_\text{corrected} := x_\text{pre} + I_{2h}^h V_{2h} (I_{2h}^h)^\trans r =  x_\text{pre} + I_{2h}^h V_{2h} (I_{2h}^h)^\trans (I - A \Lambda_1^\trans) \cdots (I - A \Lambda_\nu^\trans) b. \]
Next, apply the post-smoother to obtain
\[ x_\text{post} := (I - \Lambda_\nu A) \cdots (I - \Lambda_1 A) x_\text{corrected} + A^{-1} \Bigl[ I - (I - A \Lambda_\nu) \cdots (I - A \Lambda_1) \Bigr] b. \]
Expanding the various terms, we have
\begin{align*}
    x_\text{post} &= (I - \Lambda_\nu A) \cdots (I - \Lambda_1 A) A^{-1} \Bigl[ I - (I - A \Lambda_1^\trans) \cdots (I - A \Lambda_\nu^\trans) \Bigr] b \\
    &\quad + (I - \Lambda_\nu A) \cdots (I - \Lambda_1 A) I_{2h}^h V_{2h} (I_{2h}^h)^\trans (I - A \Lambda_1^\trans) \cdots (I - A \Lambda_\nu^\trans) b \\
    &\quad + A^{-1} \Bigl[ I - (I - A \Lambda_\nu) \cdots (I - A \Lambda_1) \Bigr] b \\
    &= \Bigl[ A^{-1} - (I - \Lambda_\nu A) \cdots (I - \Lambda_1 A) A^{-1} (I - A \Lambda_1^\trans) \cdots (I - A \Lambda_\nu^\trans) \Bigr]b  \\
    &\qquad + (I - \Lambda_\nu A) \cdots (I - \Lambda_1 A) I_{2h}^h V_{2h} (I_{2h}^h)^\trans (I - A \Lambda_1^\trans) \cdots (I - A \Lambda_\nu^\trans) b.
\end{align*}
Assuming $A$ is symmetric, the bracketed expression on the penultimate line is a symmetric matrix. (As before, if $A^{-1}$ does not exist, the expression should be interpreted symbolically in fully expanded form with the $A^{-1}$ terms appropriately removed.) Let $V$ denote the overall V-cycle operator acting on the input $b$, i.e., $b \mapsto Vb := x_\text{post}$. Assuming inductively that $V_{2h}$ is symmetric, it follows that $V$ is also symmetric. The same proof can be used in the case of a multiplicative cascading smoother by flattening $(\Lambda_{\ell,k})_{\ell,k=1}^{\nu,\chi}$ into the sequence $\Lambda_{1,1}$, $\Lambda_{1,2}, \ldots, \Lambda_{1,\chi}, \Lambda_{2,1}, \Lambda_{2,2}, \ldots$.

\section{On Preconditioned Conjugate Gradient}

A cascading smoother V-cycle that uses reverse application for the pre-smoother and forward application for the post-smoother (i.e., $\RFmode$) yields a symmetric solver, $V = V^\trans$. In the case of symmetric positive (semi)definite (SPD) problems $Ax = b$, it may be of interest, then, to use $V$ in a preconditioned Conjugate Gradient (PCG) solver. However, there is nothing in the construction of a cascading smoother that guarantees $V$ itself will be SPD.\footnote{For some Poisson test cases, using high-order ($\polydeg \geq 3$) discontinuous Galerkin methods, shallow additive cascading smoothers ($\nu \leq 2$), and $\mathsf{CS}\Jmode\RFmode$, numerical experiments confirmed that $V$ was symmetric indefinite and PCG failed. However, other problem configurations yielded SPD $V$. Whenever PCG converged, it did so with a convergence rate almost identical to preconditioned GMRES method.} To that end, one could modify the construction to build in additional constraints. For example, analysis of the V-cycle operator (see \cref{sec:Vcyclesymmetry}) suggests that $\Lambda_\ell$ should be symmetric and $I - \Lambda_\ell A$ should be a contraction. One possibility might be to compute $\Lambda_\ell$ via the standard least squares formulation, but prior to moving to the next level in the cascade, project each diagonal block of $\Lambda_\ell$ to the nearest SPD matrix satisfying the constraints, guided by a block-Jacobi method. Subsequent levels in the cascade automatically adapt to the modification because they are conditioned on the smoother error propagator of prior levels. It is unclear whether the projected cascading smoothers would degrade or perhaps even accelerate the overall multigrid solver. At the same time, the performance gains of using CG over GMRES are minimal whenever few outer iterations are needed (as is typical for highly effective multigrid preconditioners); moreover, $\FFmode$ ordering may converge faster than $\RFmode$ ordering. %
We leave this exploration open for future research.

\section{Runtime Composition on Representative Problems}
\label{sec:SM-timing}

\definecolor{barcolA}{HTML}{f5f5f5}
\definecolor{barcolB}{HTML}{427996}
\definecolor{barcolC}{HTML}{bdd7e7}
\definecolor{barcolD}{HTML}{3b3b3b}
\NewDocumentCommand{\swatch}{m}{%
  \begingroup%
  \setlength{\fboxsep}{0pt}%
  \setlength{\fboxrule}{0.5pt}%
  \fbox{\textcolor{#1}{\rule{1.2ex}{1.2ex}}}%
  \endgroup%
}

\begin{table}%
\centering%
\includegraphics[scale=0.85]{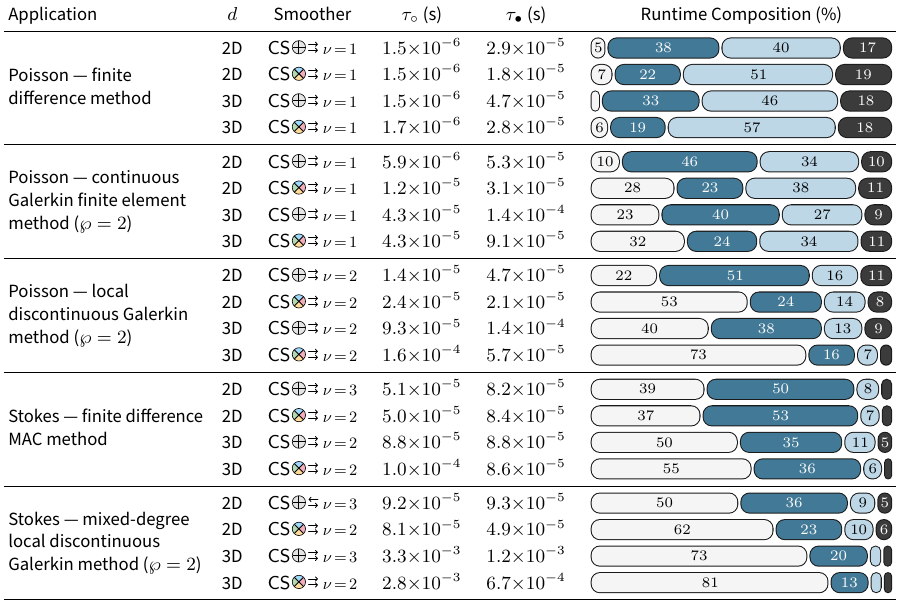}%
\caption{Empirical execution times and runtime composition for a representative set of applications. Cascading smoother build cost $\tau_\circ$ and multigrid solver cost $\tau_\bullet$ are measured in seconds per fine-grid element. See \cref{sec:SM-timing} for further details and colour mappings.}
\label{tab:o}%
\end{table}

\Cref{tab:o} provides an empirical evaluation of execution time and runtime breakdown for a representative set of applications discretised on uniform Cartesian grids of size $n \times n\,(\times\:n)$ elements. In each case, the resulting linear system is solved with a random right-hand side using multigrid V-cycle-preconditioned GMRES to a residual reduction factor of $10^{10}$. Additive and multiplicative cascading smoothers ($\mathsf{CS}\Jmode$ and $\mathsf{CS}\GSmode$) are configured with the approximately optimal ordering ($\FFmode$ or $\RFmode$) and depth ($\nu$) that minimise $\mathcal{C}_\text{rel}$. We decompose the execution time into two components: (i) the cascading smoother construction cost, accumulated across every level of the hierarchy, and (ii) the total GMRES solver time, each iteration invoking a single V-cycle. For sufficiently large $n$ and effectively stable $\eta$, both components scale approximately linearly with element count; accordingly, \cref{tab:o} reports normalised timings, with $\tau_\circ$ and $\tau_\bullet$ denoting the measured construction and solver costs per fine-grid element (in seconds), respectively.\footnote{Timing measurements were obtained on a 2017 Intel Xeon E3-1535Mv6 quad-core laptop operating at approximately 3.5\,GHz. The implementation has not been extensively tuned; for instance, the finite difference methods were emulated via a block-structured high-order DG framework.} Stacked bars illustrate the percentage breakdown of the total time $\tau_\circ + \tau_\bullet$, distinguishing the offline smoother construction ($\tau_\circ$, \swatch{barcolA}) from three online solver components: pre- and post-smoothing (\swatch{barcolB}), all other V-cycle operations such as residual restriction and correction interpolation (\swatch{barcolC}), and the GMRES driver (\swatch{barcolD}, predominantly orthogonalisation and operator evaluation $x \mapsto Ax$). As expected, the relative construction cost grows with block size, ranging from negligible for finite difference Poisson problems to appreciable for 3D Stokes problems. Depending on the application, this construction cost can be amortised across multiple solves; it can also be reduced via truncation methods, or potentially through sketching, as hypothesised in the main text.

\newcommand{\capfmt}[4]{}

\newcommand{\putfigsm}[5]{%
\par\addvspace{3ex}\noindent\begin{minipage}{\linewidth}%
\centering%
\putfig{#1}%
\captionof{figure}{\capfmt{#2}{#3}{#4}{#5}}%
\label{fig:SM-#1}%
\end{minipage}%
\par\addvspace{4ex}%
}

\section{Supplementary Numerical Experiments}

\subsection{Poisson --- Finite Difference Methods}
\label{sec:SM-Poisson-FD}

For the canonical test problem of a Poisson problem discretised by second-order finite differences on a uniform Cartesian grid, we consider two distinct approaches for the corresponding multigrid solver:
\begin{itemize}
    \item \textit{Node-centered} (results for this scheme are shown in the main article) --- In this setting, we define the interpolation operator $I_{2h}^h$ via bilinear/trilinear interpolation; e.g., for a fine-grid point located at the centre of a coarse-grid cell, its value is defined by the average of the surrounding coarse-grid points. The restriction operator $R_{h}^{2h}$ is defined by the appropriately-weighted adjoint $R_{h}^{2h} = \smash{\frac{1}{2^d} (I_{2h}^h)^\trans}$. With these choices, the Galerkin projection coarse-grid operator $R_{h}^{2h} A_h I_{2h}^h$ is not the same as the 5-point/7-point rediscretised coarse-grid operator; instead, the former is a less-sparse 9-point/27-point operator, in 2D/3D respectively. 
    
    \item \textit{Cell-centered} --- In this setting, we consider the analogue of a finite volume method and define the interpolation operator $I_{2h}^h$ via piecewise constant interpolation, i.e., the central value of a coarse-grid cell is injected into all $2^d$ fine-grid subcells. The restriction operator $R_{h}^{2h}$ is defined by the appropriately-weighted adjoint $R_{h}^{2h} = \smash{\frac{1}{2^d} (I_{2h}^h)^\trans}$. With these choices, the Galerkin projection coarse-grid operator $R_{h}^{2h} A_h I_{2h}^h$ is the same as the canonical 5-point/7-point Laplacian on the coarse grid, i.e., the Galerkin coarse-grid operator is identical to the rediscretised coarse-grid operator $A_{2h}$.    
\end{itemize}
Roughly speaking, multigrid theory predicts that piecewise constant interpolation may lead to multigrid solver breakdown because it lacks sufficient accuracy to transfer the correction from the coarse grid to the fine grid. However, this can be mitigated by applying enough smoothing after coarse-grid correction. We see this in the results for the cell-centered scheme, shown in \cref{fig:SM-T0.Nx.P1}. In particular, we observe multigrid breakdown for $\nu = 1$ in both 2D and 3D, and to a lesser extent for additive-style smoothers with $\nu = 2$ in 3D. Consequently, the overall best solvers (in terms of wall-clock time) shift from $\nu = 1$ methods in the node-centered case to $\nu \in \{2,3,4\}$ methods in the cell-centered case. On the other hand, sometimes the cell-centered scheme yields faster convergence rates than the node-centered scheme, e.g., $d = 3$, $\nu = 3$, $\mathsf{CS}\GSmode\FFmode$ yields $\eta \approx 0.7$ (node-centered) and $\eta \approx 0.5$ (cell-centered), a factor $1.4$ speedup. Besides interpolation and restriction, another contributing factor is whether the coarse-grid operator is defined via Galerkin projection or rediscretisation. In particular, a third multigrid scheme might entail using bilinear/trilinear interpolation together with Galerkin coarse-grid operators: although this approach may lead to the overall fastest convergence rates, it comes at the price of more expensive coarse-grid solves; exploring this trade-off is beyond the scope of this work.

\renewcommand{\capfmt}[4]{Multigrid solver performance for the 2D (top) and 3D (bottom) cell-centered finite difference Poisson problem considered in \cref{sec:SM-Poisson-FD}.}
\putfigsm{T0.Nx.P1}{}{}{}{}

\subsection{Poisson --- Finite Element Methods}
\label{sec:SM-Poisson-CG}

Corresponding to the continuous Galerkin, nodal finite element method considered in the main article, \cref{fig:SM-T4.Nx.P2} compiles the results for the case of $\polydeg = 1$ bilinear/trilinear elements. In this case, the blocking is trivial, for it consists of a single node per block. Note that the best-performing additive cascading smoother is faster than the fastest optimally-damped Jacobi method; in fact, in 2D the best-performing additive cascading smoother yields nearly the same $\mathcal{C}_\text{rel}$ as the best-performing multicoloured smoother. As such, it is likely the best approach in massively parallel HPC settings where the extra parallelism afforded by additive methods is advantageous, while also avoiding the need for complex graph colouring algorithms.

\renewcommand{\capfmt}[4]{Multigrid solver performance for the 2D (top) and 3D (bottom) continuous Galerkin nodal finite element method considered in \cref{sec:SM-Poisson-CG}.}
\putfigsm{T4.Nx.P2}{}{}{}{}

\subsection{Poisson --- Discontinuous Galerkin Methods}
\label{sec:SM-Poisson-DG}

Corresponding to the local discontinuous Galerkin (LDG) method considered in the main article, \cref{fig:SM-T0.Nx.P2.orbA} compiles the results for the case of bilinear/trilinear elements ($\polydeg = 1$), while \cref{fig:SM-T0.Nx.P4.orbA} compiles the results for bicubic/tricubic elements ($\polydeg = 3$). These results consider periodic boundary conditions; Dirichlet and Neumann boundary conditions yield nearly identical results. In these examples we see a stronger dependence on the ordering, i.e., all else equal, solvers that use $\RFmode$ ordering typically have slower convergence than $\FFmode$.

Meanwhile, for the adaptively refined mesh problem, \cref{fig:SM-T0.Nx.P2.orbB} compiles the results for the case of bilinear/trilinear elements ($\polydeg = 1$), while \cref{fig:SM-T0.Nx.P4.orbB} compiles the results for bicubic/tricubic elements ($\polydeg = 3$). In the bilinear/trilinear case we see that the damped Jacobi smoothers are ineffective, whereas additive cascading smoothers yield reliable convergence.

\renewcommand{\capfmt}[4]{Multigrid solver performance for the 2D (top) and 3D (bottom) uniform Cartesian grid LDG problem considered in \cref{sec:SM-Poisson-DG}, using bilinear/trilinear elements ($\polydeg = 1$).}
\putfigsm{T0.Nx.P2.orbA}{}{}{}{}

\renewcommand{\capfmt}[4]{Multigrid solver performance for the 2D (top) and 3D (bottom) uniform Cartesian grid LDG problem considered in \cref{sec:SM-Poisson-DG}, using bicubic/tricubic elements ($\polydeg = 3$).}
\putfigsm{T0.Nx.P4.orbA}{}{}{}{}

\renewcommand{\capfmt}[4]{Multigrid solver performance for the 2D (top) and 3D (bottom) nonuniform quadtree/octree mesh LDG problem considered in \cref{sec:SM-Poisson-DG}, using bilinear/trilinear elements ($\polydeg = 1$).}
\putfigsm{T0.Nx.P2.orbB}{}{}{}{}

\renewcommand{\capfmt}[4]{Multigrid solver performance for the 2D (top) and 3D (bottom) nonuniform quadtree/octree mesh LDG problem considered in \cref{sec:SM-Poisson-DG}, using bicubic/tricubic elements ($\polydeg = 3$).}
\putfigsm{T0.Nx.P4.orbB}{}{}{}{}

\subsection{Elliptic Interface --- Discontinuous Galerkin Methods}
\label{sec:SM-elliptic}

Corresponding to the multiphase elliptic interface problem considered in the main article, \cref{fig:SM-T0.Nx.P2.orbC-} through \cref{fig:SM-T0.Nx.P4.orbC+} compile the results for polynomial degrees $\polydeg \in \{1,2,3\}$ and viscosity ratios $10^{-4}$ and $10^{+4}$.

\renewcommand{\capfmt}[4]{Multigrid solver performance for the 2D (top) and 3D (bottom) high-contrast multiphase elliptic interface problem considered in \cref{sec:SM-elliptic}, viscosity ratio $10^{{#3}4}$, and #4 elements ($\polydeg = #2$).}
\newcommand{\placefig}[5]{\putfigsm{T0.Nx.P#2.orbC#4}{#1}{#3}{#4}{#5}}
\placefig{x}{2}{1}{-}{bilinear/trilinear}
\placefig{x}{2}{1}{+}{bilinear/trilinear}
\placefig{x}{3}{2}{-}{biquadratic/triquadratic}
\placefig{x}{3}{2}{+}{biquadratic/triquadratic}
\placefig{x}{4}{3}{-}{bicubic/tricubic}
\placefig{x}{4}{3}{+}{bicubic/tricubic}

\subsection{Stokes --- Taylor-Hood Finite Element Methods}
\label{sec:SM-Stokes-TH}

Corresponding to the Taylor-Hood FEM application in the main article, here we compare two normalisation choices for the Stokes prescaling mechanism (see \cref{sec:SM-prescale}): (i) $q = \infty$, the default choice, typically optimal for high-order discretisations; and (ii) $q = 2$, an intermediate alternative for medium-order discretisations, transitioning toward $q = 1$, the latter being optimal for low-order MAC staggered grids. \Cref{fig:SM-T5a.Nx.P3} and \cref{fig:SM-T5b.Nx.P3} compile the multigrid solver results for $q = \infty$ and $q = 2$, respectively; the former replicates the results from the main article. Key findings include:
\begin{itemize}
    \item Cascading smoothers with $\RFmode$ ordering are more sensitive to $q$ than those with $\FFmode$ ordering.
    \item In 3D, the optimal multicoloured cascading smoother is $\mathsf{CS}\GSmode\FFmode$, and similar convergence rates $\eta$ are obtained for both $q = \infty$ and $q = 2$. When $q = \infty$, $\mathsf{CS}\GSmode\RFmode$ performs poorly since $\eta$ deteriorates with grid refinement; choosing $q = 2$ instead yields markedly better solver performance.
    \item For multicoloured cascading smoothers in 2D, the difference between $q = \infty$ and $q = 2$ is less pronounced than in 3D, though $q = 2$ performs slightly better.
    \item In both 2D and 3D, $q = \infty$ yields the most effective additive cascading smoother, $\mathsf{CS}\Jmode\FFmode$; in 3D, its advantage over $q = 2$ is more substantial.
\end{itemize}
Overall, whether we use $q = \infty$ or $q = 2$ in the Stokes prescaling mechanism,\footnote{We also evaluated $q = 1$, but found it to be suboptimal for both $\mathsf{CS}\Jmode$ and $\mathsf{CS}\GSmode$.} cascading smoothers with $\FFmode$ ordering yield effective multigrid solvers; however, the default choice of $q = \infty$ favours $\mathsf{CS}\Jmode$, while $q = 2$ slightly benefits $\mathsf{CS}\GSmode\FFmode$ yet substantially improves $\mathsf{CS}\GSmode\RFmode$.

\renewcommand{\capfmt}[4]{Multigrid solver performance using overlapping-block smoothers for the 2D (top) and 3D (bottom) Taylor-Hood finite element Stokes problem considered in \cref{sec:SM-Stokes-TH}, using #2 in the prescaling mechanism of \cref{sec:SM-prescale}.}
\putfigsm{T5a.Nx.P3}{x}{$q = \infty$}{}{}
\putfigsm{T5b.Nx.P3}{x}{$q = 2$}{}{}

\subsection{Multiphase Stokes --- Discontinuous Galerkin Methods}
\label{sec:SM-Stokes-DG}

Corresponding to the mixed-degree LDG framework for Stokes problems and supplementing the results shown in the main article, we compile here an expanded set of test configurations, examining different kinds of boundary conditions, adaptively refined meshes, high-contrast multiphase cases, unsteady single-phase problems, and unsteady multiphase problems, as follows. We note also that cascading smoothers were tested on an \textit{equal-degree} LDG framework for Stokes problems, whereby the velocity field and pressure field have the same polynomial degree---this discretisation requires a pressure penalty stabilisation operator, but it does not alter the outcome, i.e., cascading smoothers are also highly effective in the equal-degree case.

\begin{itemize}
    \item \Cref{fig:SM-T2.Nx.P3.orbAbcP} and \cref{fig:SM-T2.Nx.P4.orbAbcP} examine single-phase, uniform-viscosity, steady-state, stress-form problems with periodic boundary conditions; the first with $\polydeg = 2$, i.e., bi/triquadratic velocity fields and bi/trilinear pressure fields; the second with $\polydeg = 3$, i.e., bi/tricubic velocity fields and bi/triquadratic pressure fields.

    \item \Cref{fig:SM-T2.Nx.P3.orbAbcN} and \cref{fig:SM-T2.Nx.P3.orbAbcD} examine single-phase, uniform-viscosity, steady-state, stress-form problems with $\polydeg = 2$; the first considers stress boundary conditions, the second Dirichlet boundary conditions. For the latter, we see an especially severe failure of block Jacobi and block Gauss-Seidel. The results for cascading smoothers show they are robust to the type of boundary condition, yielding similar results for all three kinds of boundary condition; the same is true for other polynomial degrees, though we do not include those results here.

    \item The next example tests the case of an adaptively refined mesh, identical to the geometry considered %
    in the main article for the Poisson problem. Applying stress boundary conditions (i.e., $\Gamma_D = \varnothing$ and $\Gamma_N = \partial \Omega$), \cref{fig:SM-T2.Nx.P3.orbB} compiles the multigrid solver results. In this case, we see an even stronger failure of the classical smoothers, but sufficiently deep cascading smoother-based solvers yield reliable convergence.

    \item \Cref{fig:SM-T2.Nx.P3.orbC-} and \cref{fig:SM-T2.Nx.P3.orbC+} examine high-contrast multiphase steady-state, stress-form problems with $\polydeg = 2$, viscosity ratios $10^{\pm4}$, and stress boundary conditions. Results for other polynomial degrees are similar.

    \item \Cref{fig:SM-T2.Nx.P3.orbD1} and \cref{fig:SM-T2.Nx.P3.orbD2} examine single-phase, unsteady, stress-form problems with $\polydeg = 2$ and periodic boundary conditions. In each case, $\rho = 1 / (0.1 h)$, where $h$ is the cell size of the primary mesh; this choice of $\rho$ emulates the density-weighted temporal-derivative term that would arise, e.g., from a time-stepping method with unit-order CFL. The first and second consider $\mu = 10^{-2}$ and $\mu = 10^{-4}$, representing moderately-strong and moderately-weak viscous effects, respectively. Regarding Jacobi and Gauss-Seidel, their most effective damping/overrelaxation parameters are a nontrivial function of the locally-defined grid-dependent Reynolds number scaling as $\sim \rho h / \mu$. %
    The classical smoother results of \cref{fig:SM-T2.Nx.P3.orbD1,fig:SM-T2.Nx.P3.orbD2} are tuned to the case of steady-state Stokes problems.
    As such, their solver performance is likely to improve through meticulous tuning. %
    To that end, one possibility is to optimise  %
    at $\rho \equiv 0, \mu \equiv 1$ as well as $\rho \equiv 1, \mu \equiv 0$, and then blend between these parameters depending on the specifics of the discretised problem, i.e., the effective grid-dependent Reynolds number whose strength would change across the multigrid hierarchy. Regardless, we emphasise that cascading smoothers require no such tuning and yield robust and fast convergence, across a full range of parameters.

    \item \Cref{fig:SM-T2.Nx.P3.orbEwx} and \cref{fig:SM-T2.Nx.P3.orbEgx} examine multiphase, unsteady, stress-form problems with $\polydeg = 2$ and Dirichlet boundary conditions. In these examples, density and viscosity coefficients jump by 3 to 4 orders of magnitude across an embedded interface. %
    Two scenarios are considered, a water bubble surrounded by gas, and a gas bubble surrounded by water, with $\smash{\rho_{\text{water}}} = 1$, $\smash{\mu_{\text{water}}} = 0.01$, $\smash{\rho_{\text{gas}}} = 0.001$, and $\smash{\mu_{\text{gas}}} = 0.0002$ (approximate CGS values at ambient temperature). As in the prior set of tests, the damping/overrelaxation parameters for block-Jacobi and block-Gauss-Seidel are not optimally tuned; doing so would require another parameter sweep and a nontrivial blending method between the two extremes of purely-viscous to purely-inviscid. Nevertheless, even with optimal damping, classical smoothers are expected to be ineffective since that is the case for steady-state high-contrast multiphase problems. Cascading smoothers are effective on multiphase unsteady Stokes problems; other parameter configurations and polynomial degrees yield similar outcomes.    
\end{itemize}

\renewcommand{\capfmt}[4]{Corresponding to a $\polydeg = #2$ mixed-degree LDG method, multigrid solver performance for a 2D (top) and 3D (bottom) #3.}

\putfigsm{T2.Nx.P3.orbAbcP}{x}{2}{single-phase, uniform viscosity, steady-state, stress-form Stokes problem on a uniform Cartesian grid, with periodic boundary conditions}{}
\putfigsm{T2.Nx.P4.orbAbcP}{x}{3}{single-phase, uniform viscosity, steady-state, stress-form Stokes problem on a uniform Cartesian grid, with periodic boundary conditions}{}

\putfigsm{T2.Nx.P3.orbAbcN}{x}{2}{single-phase, uniform viscosity, steady-state, stress-form Stokes problem on a uniform Cartesian grid, with stress boundary conditions, i.e., $\Gamma_D = \varnothing$ and $\Gamma_N = \partial \Omega$}{}
\putfigsm{T2.Nx.P3.orbAbcD}{x}{2}{single-phase, uniform viscosity, steady-state, stress-form Stokes problem on a uniform Cartesian grid, with Dirichlet boundary conditions, i.e., $\Gamma_N = \varnothing$ and $\Gamma_D = \partial \Omega$}{}

\putfigsm{T2.Nx.P3.orbB}{x}{2}{single-phase, uniform viscosity, steady-state, stress-form Stokes problem on a quadtree mesh (top) and octree mesh (bottom), with stress boundary conditions}{}

\putfigsm{T2.Nx.P3.orbC-}{x}{2}{high-contrast multiphase, steady-state, stress-form Stokes problem with stress boundary conditions, wherein an interior phase $\Omega_0 = (\frac14,\frac34)^d$ has viscosity $\mu_0$, the exterior phase $\Omega_1 = (0,1)^d \setminus \overline{\Omega_0}$ has viscosity $\mu_1$, and $\mu_0/\mu_1 = 10^{-4}$}{}
\putfigsm{T2.Nx.P3.orbC+}{x}{2}{high-contrast multiphase, steady-state, stress-form Stokes problem with stress boundary conditions, wherein an interior phase $\Omega_0 = (\frac14,\frac34)^d$ has viscosity $\mu_0$, the exterior phase $\Omega_1 = (0,1)^d \setminus \overline{\Omega_0}$ has viscosity $\mu_1$, and $\mu_0/\mu_1 = 10^{+4}$}{}

\putfigsm{T2.Nx.P3.orbD1}{x}{2}{single-phase, unsteady, stress-form Stokes problem with periodic boundary conditions, where $\mu = 10^{-2}$, $\rho = 1/(0.1 h)$, and $h$ is the cell size of the primary mesh}{}
\putfigsm{T2.Nx.P3.orbD2}{x}{2}{single-phase, unsteady, stress-form Stokes problem with periodic boundary conditions, where $\mu = 10^{-4}$, $\rho = 1/(0.1 h)$, and $h$ is the cell size of the primary mesh}{}

\putfigsm{T2.Nx.P3.orbEwx}{x}{2}{multiphase, unsteady, stress-form Stokes problem with Dirichlet boundary conditions, where the interior phase $\Omega_0 = (\frac14,\frac34)^d$ has $\mu_0 = 0.01$ and $\rho_0 = 1/(0.1 h)$, whereas the exterior phase $\Omega_1 = (0,1)^d \setminus \overline{\Omega_0}$ has $\mu_1 = 0.0002$ and $\rho_1 = 0.001 / (0.1 h)$, and $h$ is the cell size of the primary mesh}{}
\putfigsm{T2.Nx.P3.orbEgx}{x}{2}{multiphase, unsteady, stress-form Stokes problem with Dirichlet boundary conditions, where the interior phase $\Omega_0 = (\frac14,\frac34)^d$ has $\mu_0 = 0.0002$ and $\rho_0 = 0.001/(0.1 h)$, whereas the exterior phase $\Omega_1 = (0,1)^d \setminus \overline{\Omega_0}$ has $\mu_1 = 0.01$ and $\rho_1 = 1 / (0.1 h)$, and $h$ is the cell size of the primary mesh}{}

\fi

\end{document}